\documentclass[12pt]{elsarticle}

\usepackage{lineno}%hyperref
\usepackage{pgf,tikz}
\usepackage{mathrsfs}
\usetikzlibrary{arrows}

\usepackage{color}

\modulolinenumbers[5]

\journal{Journal of \LaTeX\ Templates}

\newtheorem{thm}{Theorem}
\newtheorem{lem}{Lemma}
\newtheorem{cor}{Corollary}
\newtheorem{pro}{Proposition}
\newdefinition{algoritmo}{Algorithm}
\newdefinition{df}{Definition}
\newdefinition{rmk}{Remark}
\newproof{pf}{Proof}
\newproof{pot}{Proof of Theorem \ref{thm2}}

\usepackage{verbatim}

\usepackage{amsmath}
\usepackage{pgf,tikz}
\usetikzlibrary{arrows}
\begin{document}

\begin{frontmatter}

\title{Cyclic structure, vertex degree and number of linear vertices in Minimal Strong Digraphs}

\author[ups,upm]{Miguel Arcos-Argudo}
\cortext[correspondingauthor]{Corresponding author.}
\ead{marcos@ups.edu.ec}

\author[upm]{Jes\'us Lacalle}
\ead{jesus.glopezdelacalle@upm.es}

\author[upm]{Luis M. Pozo-Coronado\corref{correspondingauthor}}
\ead{lm.pozo@upm.es}

\address[ups]{Grupo de Investigaci\'on en Inteligencia Artificial y
Tecnolog\'{i}as de Asistencia, Universidad Polit\'ecnica Salesiana del Ecuador}
\address[upm]{Dep. de Mat. Aplicada a las Tecnolog\'{i}as de la Informaci\'on y
las Comunicaciones, ETS de Ingenier\'{i}a de Sistemas Inform\'aticos, Universidad Polit\'ecnica de Madrid}

\begin{abstract}
Minimal Strong Digraphs (MSDs) can be regarded as a generalization of the concept of tree to
directed graphs. Their cyclic structure and some spectral properties have been studied in several
articles. In this work, we further study some properties of MSDs that have to do with bounding the
length of the longest cycle (regarding the number of linear vertices, or the maximal in-
or outdegree of vertices); studying whatever consequences from the spectral point of
view; and giving some insight about the circumstances in which an efficient algorithm to
find the longest cycle contained in an MSD can be formulated.
Among other properties, we show that the number of linear vertices contained in an MSD is
greater or equal to the maximal (resp. minimal) in- or outdegree of any vertex of the MSD
and that the maximal length of a cycle contained in an MSD is lesser or equal to $2n-m$
where $n,m$ are the order and the size of the MSD respectively; we have found a bound
for the coefficients of the characteristic polynomial of an MSD, extending the result in
\cite{cita2}, and finally, we prove that computing the longest cycle contained in an
MSD is an NP-Hard problem (Theorem \ref{complex1}).
\end{abstract}


\end{frontmatter}

%\linenumbers

\section{Introduction} \label{section1}

%%\noindent
A Minimal Strong Digraph (MSD) is a strong digraph in which the deletion of any arc yields a non strongly connected digraph. In~\cite{cita1} a compilation of properties of MSDs can be found  and in~\cite{cita2} an update of the catalog of properties is given, together with a comparative analysis between MSDs and non directed trees, where a serie of analog properties of both type of graphs is presented. In this sense, MSDs gain interest as a counterpart of trees in the context of directed graphs.

There are several other reasons to justify the interest to study MSDs. One of them is the relationship between MSDs and nearly reducible $(0,1)$-matrices (via the adjacency matrix), see, for instance, \cite{BR,HS}; and the nonnegative inverse eigenvalue problem (see \cite{TAAMP}): given real numbers $k_1, k_2,\dots, k_n$, find necessary and sufficient conditions for the existence of a nonnegative matrix $A$ of order $n$ with characteristic polynomial $x^n+k_1x^{n-1}+ k_2 x^{n-2}+\dots+ k_n$. The coefficients of the characteristic polynomial are closely related to the cycle structure of the weighted digraph with adjacency matrix $A$ by means of the Theorem of the coefficients \cite{CDS}, and the irreducible matricial realizations of the polynomial (which are identified with strongly connected digraphs \cite{BR}) can easily be reduced to the class of minimal strong digraphs. Hence, a better understanding of the cyclic structure of MSDs could lead to results on spectral theory.

Other goal for our work was trying to take advantage of the fact that minimality among SDs is a very restrictive condition. For instance, it is well-known that the size of an MSD of order $n$ is bounded by $2(n-1)$. We thought that the fact that the class of MSDs is comparatively small, together with the properties obtained in \cite{cita2}, pointing out relationships between the size of the longest cycle in an MSD and the number of linear vertices, could lead to find an algorithm of polynomial complexity to find the longest cycle in an MSD. Note that finding the
longest cycle on a SD is an NP-hard problem.

Our work plan was, thus, to further study the properties of MSDs that could give a better understanding of their cyclic structure, specially those having to do with bounding the length of the longest cycle (regarding the number of linear vertices, or the maximal in- or outdegree of vertices); studying whatever consequences from the spectral point of view; and finally trying to devise an efficient algorithm to find the longest cycle in an MSD. The first steps were accomplished, but we had to accept that the restrictions we had obtained to bound the length of cycles in an MSD were not enough to simplify the search of the longest cycle. In fact, we have proved that finding the longest cycle in an
MSD is NP-hard. Nevertheless, we think that the new properties of MSDs that we have been able to prove are interesting in themselves, insofar as they progress in the way of understanding the cyclic structure of MSDs, and hence they can lead to advances in spectral theory.

The outline of the article is as follows: In section \ref{section2} we introduce some notations and review several results on MSDs. In section \ref{section3} we study the relationship between the length of the longest cycle, the number of linear vertices, and the maximal in- or out-degree of vertices. We also state some MSDs properties, regarding chains and its contraction, that arise from the ear decomposition. In section \ref{section4}, we state a bound for the coefficients of the characteristic polynomial of an MSD, extending the results of \cite{cita2}. In section \ref{section5}, we prove that the problem of finding the longest cycle in an MSD is NP-hard. Finally, we draw some conclusions.

\section{Notation and basic properties} \label{section2}

In this paper we use some concepts and basic results about graphs that are
described below, in order to fix the notation~\cite{cita5, cita1, cita2}. %~\cite{BJG, cita1, cita2}.

Let $D=(V,A)$ be a digraph. If $(u,v)\in A$ is an arc of $D$, we say that $u$ is the initial vertex and $v$
the final vertex of the arc and we denote the arc by $uv$. We shall consider only directed paths and directed cycles. We shall denote by $n=|V|$ and by $m=|A|$ the order and the size of $D$ respectively.

In a strongly connected digraph, the \emph{indegree} $d^-(v)$ and the \emph{outdegree} $d^+(v)$ of
every vertex $v$ are greater than or equal to $1$. We shall say that $v$ is a \emph{linear vertex} if it satisfies $d^+(v)=d^-(v)=1$.

An arc $uv$ in a digraph $D$ is \emph{transitive} if there exists another $uv$-path disjoint to the arc $uv$. A digraph is
called a \emph{minimal} digraph if it has no transitive arcs.

%CONTRACCION DE UN CICLO
The \emph{contraction} of a subdigraph consists in the reduction of the subdigraph to a unique vertex $\bar{v}$. Note that the contraction of a cycle of length $q$ in an SD yields another SD. In such a process, $q-1$ vertices and $q$ arcs are eliminated. Given a cycle $C_q$, let $\bar{v}$ be the vertex corresponding to $C_q$ after contraction. We shall denote by $d^-(C_q)=d^-(\bar{v})$ (resp. $d^+(C_q)=d^+(\bar{v})$). Note that $d^+(C_q)=\sum_{v \in C_q} (d^+(v)-1)$ (and the same with $d^-(v)$).

Some basic properties concerning MDSs can be found in \cite{cita5, BJG,cita3,cita1,cita2}.

%CONTRACCION CICLOS
%MSD LINEALES
We summarize some of them: The size of an MSD digraph $D$ of order $n \geq 2$ verifies $n \leq m \leq 2(n-1)$ \cite{cita1}.
The contraction of a cycle in an MSD preserves the minimality, that is, it produces another MSD%
; hence, if we contract a strongly connected subdigraph in a minimal digraph, the resulting digraph is also minimal; each MSD of order $n\geq 2$ has at least two linear vertices.
%Construction of Rooted Tree and forest - for proof of external chain existence %
\color{black}
If $C_q$ is a cycle contained in an MSD $D$, then
the number of linear vertices of $D$ is greater than or equal to $\left\lfloor\frac{q+1}{2}\right\rfloor$.
An MSD factors into a rooted spanning tree and a forest of reversed rooted trees \cite{cita2}.
%The next result will be ratified in this work.
%Let $D$ be an MSD and $C_q$ a cycle contained in $D$. \cite{cita5}.
%(D,C_q) definition. H construction
Finally, we will use the next result.

\begin{lem}\label{lemacutpoints} \cite{cita2}
If an MSD contains a cycle $C_2$, then the vertices on the cycle are linear vertices or cut points.
\end{lem}

\color{black}

\color{black}
\section{Lower bounds of the number of linear vertices of an MSD} \label{section3}

\noindent Let $D$ be an MSD and $C_q$ a cycle contained in $D$.

In this section we show some results obtained through the analysis of the degree of the vertices, especially those with high degree.

%1 Lemma 1%%%%%%%%%%%%%%%%%%%%%%%%%%%%%%%%%%%%%%%%%%%%%%%%%%%%%%%%%%%%%%%
\color{black}
\begin{pro}\label{prop_00}%\label{lema1}
\color{black}
Let $D=(V,A)$ be an MSD, $\lambda$ the number of linear vertices of $D$ and $v\in V$ a vertex such that $v$ is contained in each cycle of $D$. Then, $\lambda \geq max(d^-(v), d^+(v))$.
\end{pro}
\color{black}
\begin{pf}
If $D$ is a cycle then $d^-(v)=d^+(v)=1$, therefore $\lambda \geq 2 > max(d^-(v),d^+(v))=1$ and the proof is completed.

Otherwise, let $C_q=v,u_1,\dots,u_{q-1},v$ be a cycle contained in $D$. By definition of MSD, each arc of $D$ is contained in at least one directed cycle of $D$, or else $D$ would not be strongly connected. Since $v$ is contained in each cycle of $D$, then each arc $wu_i$ such that $w\notin C_q$ is contained in a cycle $v,\dots,w,u_i,\dots,v$ for $1 \leq i \leq q-1$. In a similar way, each arc $u_iw$ such that $w\notin C_q$ is contained in a cycle $v,\dots,u_i,w,\dots,v$ for $1 \leq i \leq q-1$.

We shall prove that in $C_q$ there must exist at least one linear vertex.
Let us in fact \color{black} suppose, by contradiction, that $u_i \in C_q$ is not a linear vertex for $1 \leq i \leq q-1$. Hence, $d^-(u_1)=1$ or else the arc $vu_1$ would be transitive in $D$. In fact, if $d^-(u_1)>1$, since $v$ is contained in each cycle, $v$ is a vertex reached by walking in reverse direction from $u_1$ using an arc $u'_1u_1$ different from $vu_1$ (such an arc exists because of $d^{-}(u_1)>1$), and then a $vu_1$-path (not containing the arc $vu_1$) can be obtained by concatenation of a $vu'_1$-path with the arc $u'_1u_1$.

Then, $d^+(u_1)>1$, since $d^-(u_1)=1$ and we are assuming that $u_1$ is not linear. Let $u_1''\neq u_2$ be the vertex defined by the corresponding arc $u_1u''_1\in D$.

%++++
Now the following result will be proved for all $u_i$, $2\leq i\leq q-1$: $d^-(u_i)=1$ and there is an arc $u_iu_i^{\prime\prime}$ with $u_i^{\prime\prime}\neq u_{i+1}$. To show this, the following reasoning is applied iteratively for each vertex, starting from $u_2$.
First, we remark that $d^-(u_i)=1$. Otherwise, the arc $u_{i-1}u_i$ would be transitive in $D$, because an $u_{i-1}u_i$-path would exist, not containing the arc $u_{i-1}u_i$. In fact, since $v$ is contained in each cycle, $v$ is a vertex reached walking in reverse direction from $u_i$ starting with an arc $u_i^\prime u_i$ different from $u_{i-1}u_i$ (such an arc exists, since $d^-(u_i)>1$).
Also $v$ is a vertex reached walking from $u_{i-1}$ starting with the arc $u_{i-1}u_{i-1}^{\prime\prime}$.
Then an $u_{i-1}u_i$-path would be obtained by concatenation of the arc $u_{i-1}u_{i-1}^{\prime\prime}$ with the $u_{i-1}^{\prime\prime}v$-path the $vu'_{i}$-path and the arc $u'_iu_i$.

Besides, $d^+(u_i)>1$ also holds, because $d^-(u_i)=1$ and, by hypothesis, $u_i$ is not a linear vertex.
Let $u_i^{\prime\prime}\neq u_{i+1}$ be the vertex defined by the arc $u_iu_i^{\prime\prime}$ belonging to $D$.

Finally, let us show that the arc $u_{q-1}v$ is transitive. In fact, since $v$ is contained in each cycle, $v$ is a vertex reached walking from $u_{q-1}$, starting with the arc $u_{q-1}u_{q-1}^{\prime\prime}$. The $u_{q-1}v$-path obtained by concatenation of the arc $u_{q-1}u_{q-1}^{\prime\prime}$ with the $u_{q-1}^{\prime\prime}v$-path proves that $u_{q-1}v$ is transitive. This fact contradicts the minimality of $D$.

We have still to prove that the linear vertices reached for each outgoing (resp. incoming) arc from (resp. to) $v$ are all different. Let $vu_1$ and $vu'_1$ be two arcs in $D$.
From $vu_1$, as we have seen, we can construct a path $v,u_{1},\dots,u_{k}$ such that $d^{-}(u_{i})=1$ for $1 \leq i \leq k$ and $d^{+}(u_{i})>1$ for $1 \leq i \leq k-1$ and $u_{k}$ is linear (note that $k$ can be 1, but it must exist) as we have proved previously. Now, in a similar way, we construct a path $v,u'_{1},\dots,u'_{l}$ such that $d^{-}(u'_{j})=1$ for $1 \leq i \leq l$ and $d^{+}(u'_{j})>1$ for $1 \leq i \leq l-1$ and $u'_{l}$ is linear.

The paths $v,u_{1},\dots,u_{k}$ and $v,u'_{1},\dots,u'_{l}$ cannot rejoin after they leave $v$, since all the indegrees of their vertices are $1$. Hence, $u_k \neq u'_l$. The proof is completed.
\color{black}

\end{pf}

%1.1 a Lemma 3%%%%%%%%%%%%%%%%%%%%%%%%%%%%%%%%%%%%%%%%%%%%%%%%%%%%%%%%%%%%%%%
\color{black}
\begin{pro}\label{prop_0}%\label{lema1_1}
\color{black}
Let $D=(V,A)$ be an MSD of order $n \geq 2$, $v \in V$ a vertex of $D$ and $\lambda$ the number of linear vertices contained in $D$. Then, $\lambda \geq max(d^-(v), d^+(v))$.
\end{pro}
\color{black}
\begin{pf}
If $D=C_n$ then $d^-(v)=d^+(v)=1$, therefore $\lambda=q\geq2$ and the proof is completed.

Otherwise, we obtain an MSD $D'$ from $D$ by the contraction of all cycles that don't contain the vertex $v$. Note that $v$ is a vertex contained in each cycle of $D'$. Then, by \color{black} Proposition \ref{prop_00} \color{black} $\lambda_{D'} \geq max(d^-(v),d^+(v))$ where $\lambda_{D'}$ is the number of linear vertices of $D'$. Note also that $v$ preserves in $D'$ all its incident arcs. \color{black}
Next we expand the cycles contracted previously. In this process the linear vertices are maintained. Indeed, if we expand a linear vertex corresponding to a cycle of length greater than 2 this fact is obvious. And, if we expand one corresponding to a cycle of length 2 the result follows from Lemma \ref{lemacutpoints}.
The proof is completed.
\end{pf}
\color{black}

%1.3 a Lemma 4%%%%%%%%%%%%%%%%%%%%%%%%%%%%%%%%%%%%%%%%%%%%%%%%%%%%%%%%%%%%%%%
%Contraction of cycles (ok)

\begin{cor}\label{lema1_3}
Let $D=(V,A)$ be an MSD, $C_q$ a cycle contained in $D$ and $\mu$ the number of linear vertices contained in $D$ but not contained in $C_q$. Then, $\mu \geq max(d^-(C_q),d^+(C_q))$.
\end{cor}
\begin{pf}
If $D=C_q$ then $\mu = d^-(C_q) = d^+(C_q) = 0$ and the proof is completed.

%Otherwise, we obtain an MSD $D'$ from $D$ by first contracting $C_q$ in a unique vertex $v'$ and then contracting all cycles that not contain $v'$. Note that $v'$ is a vertex contained in each cycle of $D'$. Then, by \color{black} Proposition \ref{prop_00}\color{black}, $\mu' \geq max(d^-(v'),d^+(v'))$ where $\mu'$ is the number of linear vertices of $D'$.
%%Now, we obtain an MSD $D'$ from $\bar{D}$ by the expansion of the cycle $C_q$ (contracted in $\bar{v}$).
%Note that $d^-(C_q)=d^-(v')$ and $d^+(C_q)=d^+(v')$ because $v'$ preserves all incident arcs on vertices of $C_q$, except the arcs of $C_q$. Then $\mu' \geq max(d^-(C_q),d^+(C_q))$.
%%Now, in $D'$ we suppress all arcs of $C_q$ to obtain a digraph $D''$ that is a digraph and a Hasse diagram associated to $(D',C_q)$.
%Finally, in the same way as in the proof of \color{black} Proposition \ref{prop_0} \color{black}, each linear vertex $u \in D'$ %XXX
%is a linear vertex in $D$ or is a cycle that contains a linear vertex in $D$. The proof is thus completed.
Otherwise, we obtain an MSD $D'$ from $D$ by contracting $C_q$ in a unique vertex $v'$. Note that the number of linear vertices of $D'$ is precisely $\mu$. Application of Proposition~\ref{prop_0} then implies that $\mu\ge max(d^-(v'),d^+(v'))=max(d^-(C),d^+(C))$ and we are finished.
\end{pf}

%1.2 a Corollary 1%%%%%%%%%%%%%%%%%%%%%%%%%%%%%%%%%%%%%%%%%%%%%%%%%%%%%%%%%%%%%%%
As we mentioned in \ref{section2}, if there is a cycle $C_q \in D$, the number of linear vertices of $D$ is greater than or equal to $\left\lfloor\frac{q+1}{2}\right\rfloor$, see \cite{cita5}. \color{black} We ratify this result with a new, shorter proof, by using previous properties.
\begin{cor}\label{cor1_2}
Let $D=(V,A)$ be an MSD of order $n \geq 2$, $C_q$ a cycle contained in $D$ and $\lambda$ the number of linear vertices contained in $D$. Then, $\lambda \geq \left\lfloor\frac{q+1}{2}\right\rfloor$.
\end{cor}
\begin{pf}
Let $\nu$ be the number of linear vertices contained in $C_q$ and $\mu$ the rest of linear vertices of $D$. Then $\lambda=\mu + \nu$, and we know by Corollary \ref{lema1_3} that $\mu \geq max(d^+(C_q),d^-(C_q))$. Since $d^{+}(C_q)+d^-(C_q) \geq q - \nu$ we have that

$$\mu \geq max(d^+(C_q),d^-(C_q)) \geq \left\lceil \frac{q - \nu}{2} \right\rceil,$$

and then

$$\lambda=\mu + \nu \geq \left\lceil \frac{q-\nu}{2} \right\rceil + \nu = \left\lceil \frac{q+\nu}{2} \right\rceil \geq \left\lceil \frac{q}{2} \right\rceil = \left\lfloor \frac{q+1}{2} \right\rfloor.$$
The proof is completed.
\end{pf}

%\color{blue}

%Corollary acotacion inversa de la longitud del ciclo mas largo%%%%%%%%%%
As a consequence of Corollary \ref{cor1_2} we obtain an upper bound for the maximum length of a cycle contained in an MSD.
\begin{cor}\label{cor_bound_inversa}
Let $D=(V,A)$ be an MSD of order $n \geq 2$, $C_l$ a cycle with maximal length $l$ contained in $D$ and $\lambda$ the number of linear vertices contained in $D$. Then, $l \leq 2\lambda$.
\end{cor}
\begin{pf}
By Corollary \ref{cor1_2} we know that
$$\lambda \geq \left\lfloor\frac{l+1}{2}\right\rfloor,$$
then
$$l \leq 2\lambda.$$
The proof is completed.
\end{pf}

%**************************************************************************************

Since every vertex contained in an MSD must be contained in at least one directed cycle, we can obtain two different bounds for the number of linear vertices, one from the vertex degree and one from the cycle length. The next result somehow combines the two aforementioned bounds.

\begin{cor}\label{cor_lamda_maxgrado}
Let $D=(V,A)$ be an MSD of order $n \geq 2$, $C_q$ a directed cycle of length $q$ contained in $D$, $u \in C_q$ a vertex of $D$ and $\lambda$ the number of linear vertices contained in $D$. Then,
$$\lambda \geq \left\lfloor\frac{q+d(u)}{2}\right\rfloor-1.$$
\end{cor}
\begin{pf}
As we did in the proof of Corollary~\ref{cor1_2}, we call $\nu$ the number of linear vertices contained in $C_q$ and $\mu$ the rest of linear vertices of $D$. $\nu$ attains its minimum if the in- and outdegree of $C_q$ are maximally spread out among the vertices of $C_q$. But now we now that part of the in- and outdegree of $C_q$ is condensed in the vertex $u$. We thus obtain the inequality
\begin{multline*}
\nu + (d^+(C_q)-(d^+(u)-1))+(d^-(C_q)-(d^-(u)-1))+1 \geq q
\\ \Rightarrow d^+(C_q)+d^-(C_q) \geq q -\nu + d(u) - 3.
\end{multline*}
Combining it with Corollary~\label{lema1_3} ($\mu \geq max(d^+(C_q),d^-(C_q))$),
we obtain
\[
\mu \geq max(d^+(C_q),d^-(C_q)) \geq \left\lceil \frac{q -\nu + d(u) -3 }{2} \right\rceil,
\]
and finally
\begin{multline*}
\lambda=\mu + \nu \geq \left\lceil \frac{q -\nu + d(u) -3}{2} \right\rceil + \nu = \left\lceil \frac{q+\nu+ d(u) - 1}{2} \right\rceil -1
  \\\geq \left\lceil \frac{q+d(u)-1}{2} \right\rceil-1 = \left\lfloor \frac{q+d(u)}{2} \right\rfloor-1.
\end{multline*}
The proof is completed.
%If $D=C_q$ then $\lambda = q$ and the proof is completed.

%Otherwise, by Proposition \ref{prop_0} we know that $\lambda \geq max(d^-(u),d^+(u))$ and $u$ is a not linear vertex. Furthermore, by Corollary \ref{lema1_3} we know that $\mu \geq max(d^-(C_q),d^+(C_q))$ where $\mu$ is the number of linear vertices contained in $D$ but not contained in $C_q$. Now, we can suppose that $d^+(u) > d^-(u)$ and then the number of output arcs of $C_q$ from $u$ is $d^+(u)-1$. This implies that there are at least $d^+(u)-1$ linear vertices that are contained in $D$ but are not contained in $C_q$, and so
%
%$$\lambda \geq \left\lfloor\frac{q+1}{2}\right\rfloor +d^+(u)-1.$$
%
%If $d^-(u) > d^+(u)$ the proof is analogous. The proof is completed.
\end{pf}

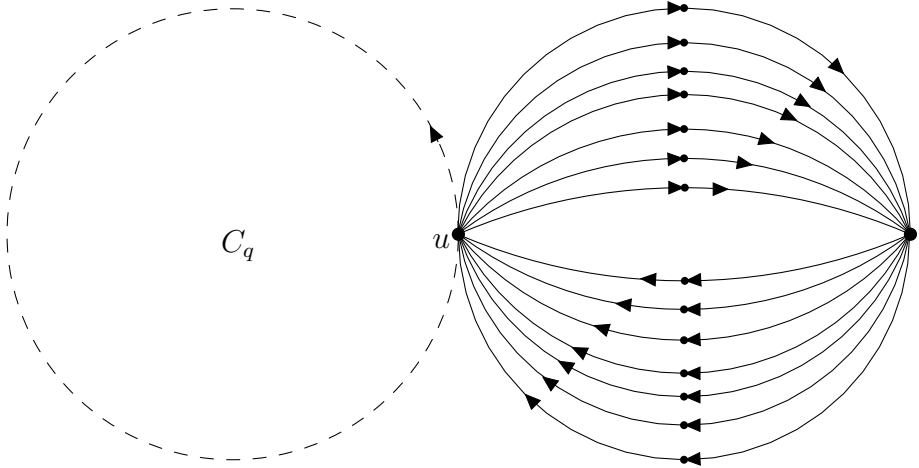
\begin{figure}[hp]
\centering
\begin{tikzpicture}[line cap=round,line join=round,>=triangle 45,x=1.0cm,y=1.0cm]
\clip(4.81,-8.17) rectangle (17.25,-1.78);
\draw [dash pattern=on 5pt off 5pt] (8,-5) circle (3cm);
\draw (7.7,-4.8) node[anchor=north west] {$C_q$};
\draw (10.51,-4.85) node[anchor=north west] {$u$};
\draw [shift={(14,2)}] plot[domain=4.31:5.12,variable=\t]({1*7.62*cos(\t r)+0*7.62*sin(\t r)},{0*7.62*cos(\t r)+1*7.62*sin(\t r)});
\draw [shift={(14,-1)}] plot[domain=4.07:5.36,variable=\t]({1*5*cos(\t r)+0*5*sin(\t r)},{0*5*cos(\t r)+1*5*sin(\t r)});
\draw [shift={(14,-2.5)}] plot[domain=3.84:5.59,variable=\t]({1*3.91*cos(\t r)+0*3.91*sin(\t r)},{0*3.91*cos(\t r)+1*3.91*sin(\t r)});
\draw [shift={(14,-3.5)}] plot[domain=3.61:5.82,variable=\t]({1*3.35*cos(\t r)+0*3.35*sin(\t r)},{0*3.35*cos(\t r)+1*3.35*sin(\t r)});
\draw [shift={(14,-4)}] plot[domain=3.46:5.96,variable=\t]({1*3.16*cos(\t r)+0*3.16*sin(\t r)},{0*3.16*cos(\t r)+1*3.16*sin(\t r)});
\draw [shift={(14,-4.5)}] plot[domain=3.31:6.12,variable=\t]({1*3.04*cos(\t r)+0*3.04*sin(\t r)},{0*3.04*cos(\t r)+1*3.04*sin(\t r)});
\draw [shift={(14,-5)}] plot[domain=-3.14:0,variable=\t]({1*3*cos(\t r)+0*3*sin(\t r)},{0*3*cos(\t r)+1*3*sin(\t r)});
\draw [shift={(14,-5)}] plot[domain=0:3.1416,variable=\t]({1*3.01*cos(\t r)+0*3.01*sin(\t r)},{0*3.01*cos(\t r)+1*3.01*sin(\t r)});
\draw [shift={(14,-5.5)}] plot[domain=0.16:2.98,variable=\t]({1*3.05*cos(\t r)+0*3.05*sin(\t r)},{0*3.05*cos(\t r)+1*3.05*sin(\t r)});
\draw [shift={(14,-6)}] plot[domain=0.32:2.82,variable=\t]({1*3.17*cos(\t r)+0*3.17*sin(\t r)},{0*3.17*cos(\t r)+1*3.17*sin(\t r)});
\draw [shift={(14,-6.5)}] plot[domain=0.46:2.68,variable=\t]({1*3.36*cos(\t r)+0*3.36*sin(\t r)},{0*3.36*cos(\t r)+1*3.36*sin(\t r)});
\draw [shift={(14,-7.54)}] plot[domain=0.7:2.44,variable=\t]({1*3.94*cos(\t r)+0*3.94*sin(\t r)},{0*3.94*cos(\t r)+1*3.94*sin(\t r)});
\draw [shift={(14,-9)}] plot[domain=0.93:2.21,variable=\t]({1*5.01*cos(\t r)+0*5.01*sin(\t r)},{0*5.01*cos(\t r)+1*5.01*sin(\t r)});
\draw [shift={(14,-12)}] plot[domain=1.16:1.98,variable=\t]({1*7.62*cos(\t r)+0*7.62*sin(\t r)},{0*7.62*cos(\t r)+1*7.62*sin(\t r)});
\draw [->] (13.83,-1.99) -- (14,-1.99);
\draw [->] (13.85,-2.45) -- (14,-2.45);
\draw [->] (13.84,-2.84) -- (14,-2.83);
\draw [->] (13.83,-3.14) -- (14,-3.14);
\draw [->] (13.85,-3.6) -- (14,-3.6);
\draw [->] (13.83,-4) -- (14,-3.99);
\draw [->] (13.84,-4.38) -- (14.01,-4.38);
\draw [->] (14.17,-5.61) -- (14.01,-5.62);
\draw [->] (14.17,-6) -- (14,-6);
\draw [->] (14.18,-6.4) -- (14,-6.41);
\draw [->] (14.18,-6.85) -- (14,-6.85);
\draw [->] (14.18,-7.16) -- (14,-7.16);
\draw [->] (14.18,-7.54) -- (14,-7.54);
\draw [->] (14.19,-7.99) -- (14,-8);
\draw [->] (16.03,-2.77) -- (16.13,-2.87);
\draw [->] (15.76,-3) -- (15.89,-3.11);
\draw [->] (15.55,-3.23) -- (15.69,-3.31);
\draw [->] (15.32,-3.4) -- (15.51,-3.49);
\draw [->] (15.05,-3.74) -- (15.21,-3.79);
\draw [->] (14.74,-4.05) -- (14.92,-4.08);
\draw [->] (14.49,-4.4) -- (14.6,-4.4);
\draw [->] (13.56,-5.6) -- (13.41,-5.59);
\draw [->] (13.26,-5.94) -- (13.08,-5.92);
\draw [->] (12.91,-6.25) -- (12.79,-6.21);
\draw [->] (12.67,-6.58) -- (12.5,-6.5);
\draw [->] (12.44,-6.75) -- (12.32,-6.68);
\draw [->] (12.23,-6.98) -- (12.11,-6.89);
\draw [->] (11.99,-7.23) -- (11.88,-7.12);
\draw [->] (10.73,-3.77) -- (10.62,-3.55);
\begin{scriptsize}
\fill [color=black] (11,-5) circle (2.5pt);
\fill [color=black] (17.01,-5) circle (2.5pt);
\fill [color=black] (14.01,-5.62) circle (1.5pt);
\fill [color=black] (14,-6) circle (1.5pt);
\fill [color=black] (14,-6.41) circle (1.5pt);
\fill [color=black] (14,-6.85) circle (1.5pt);
\fill [color=black] (14,-7.16) circle (1.5pt);
\fill [color=black] (14,-8) circle (1.5pt);
\fill [color=black] (14,-7.54) circle (1.5pt);
\fill [color=black] (14.01,-4.38) circle (1.5pt);
\fill [color=black] (14,-3.99) circle (1.5pt);
\fill [color=black] (14,-3.6) circle (1.5pt);
\fill [color=black] (14,-3.14) circle (1.5pt);
\fill [color=black] (14,-2.83) circle (1.5pt);
\fill [color=black] (14,-2.45) circle (1.5pt);
\fill [color=black] (14,-1.99) circle (1.5pt);
\fill [color=black] (11,-5.02) circle (1.5pt);
\end{scriptsize}
\end{tikzpicture}
\caption{Example 1 for Corollary \ref{cor_lamda_maxgrado}.}
\label{fig_ejemplo1}
\end{figure}

\begin{figure}[hp]
\centering
\begin{tikzpicture}[line cap=round,line join=round,>=triangle 45,x=1.0cm,y=1.0cm]
\clip(4.68,-15.53) rectangle (17.47,-8.42);
\draw [dash pattern=on 6pt off 6pt] (8,-5) circle (3cm);
\draw (7.61,-4.43) node[anchor=north west] {$C_q$};
\draw (10.52,-4.4) node[anchor=north west] {$v$};
\draw [shift={(14,2)}] plot[domain=4.31:5.12,variable=\t]({1*7.62*cos(\t r)+0*7.62*sin(\t r)},{0*7.62*cos(\t r)+1*7.62*sin(\t r)});
\draw [shift={(14,-1)}] plot[domain=4.07:5.36,variable=\t]({1*5*cos(\t r)+0*5*sin(\t r)},{0*5*cos(\t r)+1*5*sin(\t r)});
\draw [shift={(14,-2.5)}] plot[domain=3.84:5.59,variable=\t]({1*3.91*cos(\t r)+0*3.91*sin(\t r)},{0*3.91*cos(\t r)+1*3.91*sin(\t r)});
\draw [shift={(14,-3.5)}] plot[domain=3.61:5.82,variable=\t]({1*3.35*cos(\t r)+0*3.35*sin(\t r)},{0*3.35*cos(\t r)+1*3.35*sin(\t r)});
\draw [shift={(14,-4)}] plot[domain=3.46:5.96,variable=\t]({1*3.16*cos(\t r)+0*3.16*sin(\t r)},{0*3.16*cos(\t r)+1*3.16*sin(\t r)});
\draw [shift={(14,-4.5)}] plot[domain=3.31:6.12,variable=\t]({1*3.04*cos(\t r)+0*3.04*sin(\t r)},{0*3.04*cos(\t r)+1*3.04*sin(\t r)});
\draw [shift={(14,-5)}] plot[domain=-3.14:0,variable=\t]({1*3*cos(\t r)+0*3*sin(\t r)},{0*3*cos(\t r)+1*3*sin(\t r)});
\draw [shift={(14,-5)}] plot[domain=0:3.1416,variable=\t]({1*3.01*cos(\t r)+0*3.01*sin(\t r)},{0*3.01*cos(\t r)+1*3.01*sin(\t r)});
\draw [shift={(14,-5.5)}] plot[domain=0.16:2.98,variable=\t]({1*3.05*cos(\t r)+0*3.05*sin(\t r)},{0*3.05*cos(\t r)+1*3.05*sin(\t r)});
\draw [shift={(14,-6)}] plot[domain=0.32:2.82,variable=\t]({1*3.17*cos(\t r)+0*3.17*sin(\t r)},{0*3.17*cos(\t r)+1*3.17*sin(\t r)});
\draw [shift={(14,-6.5)}] plot[domain=0.46:2.68,variable=\t]({1*3.36*cos(\t r)+0*3.36*sin(\t r)},{0*3.36*cos(\t r)+1*3.36*sin(\t r)});
\draw [shift={(14,-7.54)}] plot[domain=0.7:2.44,variable=\t]({1*3.94*cos(\t r)+0*3.94*sin(\t r)},{0*3.94*cos(\t r)+1*3.94*sin(\t r)});
\draw [shift={(14,-9)}] plot[domain=0.93:2.21,variable=\t]({1*5.01*cos(\t r)+0*5.01*sin(\t r)},{0*5.01*cos(\t r)+1*5.01*sin(\t r)});
\draw [shift={(14,-12)}] plot[domain=1.16:1.98,variable=\t]({1*7.62*cos(\t r)+0*7.62*sin(\t r)},{0*7.62*cos(\t r)+1*7.62*sin(\t r)});
\draw [->] (13.83,-1.99) -- (14,-1.99);
\draw [->] (13.85,-2.45) -- (14,-2.45);
\draw [->] (13.84,-2.84) -- (14,-2.83);
\draw [->] (13.83,-3.14) -- (14,-3.14);
\draw [->] (13.85,-3.6) -- (14,-3.6);
\draw [->] (13.83,-4) -- (14,-3.99);
\draw [->] (13.84,-4.38) -- (14.01,-4.38);
\draw [->] (14.17,-5.61) -- (14.01,-5.62);
\draw [->] (14.17,-6) -- (14,-6);
\draw [->] (14.18,-6.4) -- (14,-6.41);
\draw [->] (14.18,-6.85) -- (14,-6.85);
\draw [->] (14.18,-7.16) -- (14,-7.16);
\draw [->] (14.18,-7.54) -- (14,-7.54);
\draw [->] (14.19,-7.99) -- (14,-8);
\draw [->] (16.03,-2.77) -- (16.13,-2.87);
\draw [->] (15.76,-3) -- (15.89,-3.11);
\draw [->] (15.55,-3.23) -- (15.69,-3.31);
\draw [->] (15.32,-3.4) -- (15.51,-3.49);
\draw [->] (15.05,-3.74) -- (15.21,-3.79);
\draw [->] (14.74,-4.05) -- (14.92,-4.08);
\draw [->] (14.49,-4.4) -- (14.6,-4.4);
\draw [->] (13.56,-5.6) -- (13.41,-5.59);
\draw [->] (13.26,-5.94) -- (13.08,-5.92);
\draw [->] (12.91,-6.25) -- (12.79,-6.21);
\draw [->] (12.67,-6.58) -- (12.5,-6.5);
\draw [->] (12.44,-6.75) -- (12.32,-6.68);
\draw [->] (12.23,-6.98) -- (12.11,-6.89);
\draw [->] (11.99,-7.23) -- (11.88,-7.12);
\draw [->] (10.73,-3.77) -- (10.62,-3.55);
\draw [dash pattern=on 6pt off 6pt] (8,-12) circle (3cm);
\draw (7.7,-11.8) node[anchor=north west] {$C_q$};
\draw (10.4,-11.8) node[anchor=north west] {$u$};
\draw [shift={(14,-5)}] plot[domain=4.31:5.12,variable=\t]({1*7.62*cos(\t r)+0*7.62*sin(\t r)},{0*7.62*cos(\t r)+1*7.62*sin(\t r)});
\draw [shift={(14,-8)}] plot[domain=4.07:5.36,variable=\t]({1*5*cos(\t r)+0*5*sin(\t r)},{0*5*cos(\t r)+1*5*sin(\t r)});
\draw [shift={(14,-9.5)}] plot[domain=3.84:5.59,variable=\t]({1*3.91*cos(\t r)+0*3.91*sin(\t r)},{0*3.91*cos(\t r)+1*3.91*sin(\t r)});
\draw [shift={(14,-10.5)}] plot[domain=3.61:5.82,variable=\t]({1*3.35*cos(\t r)+0*3.35*sin(\t r)},{0*3.35*cos(\t r)+1*3.35*sin(\t r)});
\draw [shift={(14,-11)}] plot[domain=3.46:5.96,variable=\t]({1*3.16*cos(\t r)+0*3.16*sin(\t r)},{0*3.16*cos(\t r)+1*3.16*sin(\t r)});
\draw [shift={(14,-11.5)}] plot[domain=3.31:6.12,variable=\t]({1*3.04*cos(\t r)+0*3.04*sin(\t r)},{0*3.04*cos(\t r)+1*3.04*sin(\t r)});
\draw [shift={(14,-12)}] plot[domain=-3.14:0,variable=\t]({1*3*cos(\t r)+0*3*sin(\t r)},{0*3*cos(\t r)+1*3*sin(\t r)});
\draw [shift={(14,-12)}] plot[domain=0:3.1416,variable=\t]({1*3.01*cos(\t r)+0*3.01*sin(\t r)},{0*3.01*cos(\t r)+1*3.01*sin(\t r)});
\draw [shift={(14,-12.5)}] plot[domain=0.16:2.98,variable=\t]({1*3.05*cos(\t r)+0*3.05*sin(\t r)},{0*3.05*cos(\t r)+1*3.05*sin(\t r)});
\draw [shift={(14,-13.5)}] plot[domain=0.46:2.68,variable=\t]({1*3.36*cos(\t r)+0*3.36*sin(\t r)},{0*3.36*cos(\t r)+1*3.36*sin(\t r)});
\draw [shift={(14,-14.54)}] plot[domain=0.7:2.44,variable=\t]({1*3.94*cos(\t r)+0*3.94*sin(\t r)},{0*3.94*cos(\t r)+1*3.94*sin(\t r)});
\draw [shift={(14,-16)}] plot[domain=0.93:2.21,variable=\t]({1*5.01*cos(\t r)+0*5.01*sin(\t r)},{0*5.01*cos(\t r)+1*5.01*sin(\t r)});
\draw [shift={(14,-19)}] plot[domain=1.16:1.98,variable=\t]({1*7.62*cos(\t r)+0*7.62*sin(\t r)},{0*7.62*cos(\t r)+1*7.62*sin(\t r)});
\draw [->] (13.83,-8.99) -- (14,-8.99);
\draw [->] (13.85,-9.45) -- (14,-9.45);
\draw [->] (13.83,-10.14) -- (14,-10.14);
\draw [->] (13.85,-10.6) -- (14,-10.6);
\draw [->] (13.83,-11) -- (14,-10.99);
\draw [->] (13.84,-11.38) -- (14.01,-11.38);
\draw [->] (16.03,-9.77) -- (16.13,-9.87);
\draw [->] (15.76,-10) -- (15.89,-10.11);
\draw [->] (15.32,-10.4) -- (15.51,-10.49);
\draw [->] (15.05,-10.74) -- (15.21,-10.79);
\draw [->] (14.74,-11.05) -- (14.92,-11.08);
\draw [->] (14.49,-11.4) -- (14.6,-11.4);
\draw [->] (10.73,-10.77) -- (10.62,-10.55);
\draw [->] (13.95,-15) -- (14,-15);
\draw [->] (13.94,-14.54) -- (14,-14.54);
\draw [->] (13.97,-14.16) -- (14,-14.16);
\draw [->] (13.97,-13.85) -- (14,-13.85);
\draw [->] (13.97,-13.4) -- (14,-13.41);
\draw [->] (13.97,-13) -- (14,-13);
\draw [->] (13.97,-12.62) -- (14,-12.62);
\draw [->] (14.33,-12.61) -- (14.47,-12.6);
\draw [->] (14.57,-12.97) -- (14.76,-12.94);
\draw [->] (14.82,-13.32) -- (15.03,-13.27);
\draw [->] (15.08,-13.67) -- (15.3,-13.59);
\draw [->] (15.28,-13.89) -- (15.47,-13.8);
\draw [->] (15.51,-14.14) -- (15.67,-14.04);
\draw [->] (15.74,-14.44) -- (15.9,-14.32);
\draw [shift={(14,-12.87)}] plot[domain=0.28:2.86,variable=\t]({1*3.13*cos(\t r)+0*3.13*sin(\t r)},{0*3.13*cos(\t r)+1*3.13*sin(\t r)});
\draw [->] (13.86,-9.74) -- (14,-9.74);
\draw [->] (15.58,-10.16) -- (15.74,-10.26);
\draw (17.01,-12)-- (11,-12);
\draw [->] (14.07,-12) -- (13.73,-12);
\begin{scriptsize}
\fill [color=black] (11,-5) circle (2.5pt);
\fill [color=black] (17.01,-5) circle (2.5pt);
\fill [color=black] (14.01,-5.62) circle (1.5pt);
\fill [color=black] (14,-6) circle (1.5pt);
\fill [color=black] (14,-6.41) circle (1.5pt);
\fill [color=black] (14,-6.85) circle (1.5pt);
\fill [color=black] (14,-7.16) circle (1.5pt);
\fill [color=black] (14,-8) circle (1.5pt);
\fill [color=black] (14,-7.54) circle (1.5pt);
\fill [color=black] (14.01,-4.38) circle (1.5pt);
\fill [color=black] (14,-3.99) circle (1.5pt);
\fill [color=black] (14,-3.6) circle (1.5pt);
\fill [color=black] (14,-3.14) circle (1.5pt);
\fill [color=black] (14,-2.83) circle (1.5pt);
\fill [color=black] (14,-2.45) circle (1.5pt);
\fill [color=black] (14,-1.99) circle (1.5pt);
\fill [color=black] (11,-5.02) circle (1.5pt);
\fill [color=black] (11,-12) circle (2.5pt);
\fill [color=black] (17.01,-12) circle (2.5pt);
\fill [color=black] (14.01,-11.38) circle (1.5pt);
\fill [color=black] (14,-10.99) circle (1.5pt);
\fill [color=black] (14,-10.6) circle (1.5pt);
\fill [color=black] (14,-10.14) circle (1.5pt);
\fill [color=black] (14,-9.45) circle (1.5pt);
\fill [color=black] (14,-8.99) circle (1.5pt);
\fill [color=black] (11,-12) circle (1.5pt);
\fill [color=black] (14,-12.62) circle (1.5pt);
\fill [color=black] (14,-13) circle (1.5pt);
\fill [color=black] (14,-13.41) circle (1.5pt);
\fill [color=black] (14,-13.85) circle (1.5pt);
\fill [color=black] (14,-14.16) circle (1.5pt);
\fill [color=black] (14,-14.54) circle (1.5pt);
\fill [color=black] (14,-15) circle (1.5pt);
\fill [color=black] (14,-9.74) circle (1.5pt);
\end{scriptsize}
\end{tikzpicture}
\caption{Example 2 for Corollary \ref{cor_lamda_maxgrado}.}
\label{fig_ejemplo2}
\end{figure}
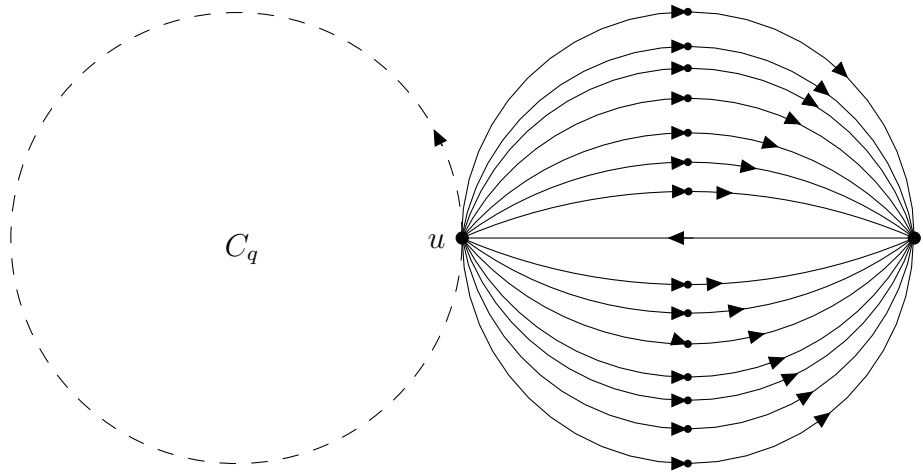

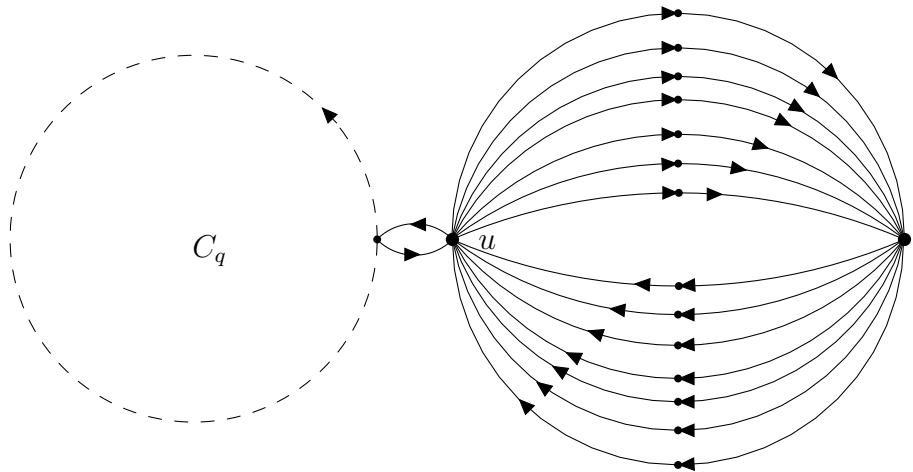
\begin{figure}[hp]
\centering
\begin{tikzpicture}[line cap=round,line join=round,>=triangle 45,x=1.0cm,y=1.0cm]
\clip(4.85,-8.3) rectangle (17.24,-1.79);
\draw [dash pattern=on 5pt off 5pt] (8,3) circle (3cm);
\draw [shift={(18,4.5)}] plot[domain=2.78:3.5,variable=\t]({1*4.27*cos(\t r)+0*4.27*sin(\t r)},{0*4.27*cos(\t r)+1*4.27*sin(\t r)});
\draw [shift={(10,4.5)}] plot[domain=-0.36:0.36,variable=\t]({1*4.27*cos(\t r)+0*4.27*sin(\t r)},{0*4.27*cos(\t r)+1*4.27*sin(\t r)});
\draw [shift={(15.5,7)}] plot[domain=4.35:5.07,variable=\t]({1*4.27*cos(\t r)+0*4.27*sin(\t r)},{0*4.27*cos(\t r)+1*4.27*sin(\t r)});
\draw [shift={(15.5,-1)}] plot[domain=1.21:1.93,variable=\t]({1*4.27*cos(\t r)+0*4.27*sin(\t r)},{0*4.27*cos(\t r)+1*4.27*sin(\t r)});
\draw [shift={(10,1.5)}] plot[domain=-0.36:0.36,variable=\t]({1*4.27*cos(\t r)+0*4.27*sin(\t r)},{0*4.27*cos(\t r)+1*4.27*sin(\t r)});
\draw [shift={(17.31,1.79)}] plot[domain=1.91:2.79,variable=\t]({1*3.52*cos(\t r)+0*3.52*sin(\t r)},{0*3.52*cos(\t r)+1*3.52*sin(\t r)});
\draw [shift={(12.62,6.54)}] plot[domain=5.08:5.9,variable=\t]({1*3.8*cos(\t r)+0*3.8*sin(\t r)},{0*3.8*cos(\t r)+1*3.8*sin(\t r)});
\draw [shift={(12.74,-0.38)}] plot[domain=0.35:1.21,variable=\t]({1*3.6*cos(\t r)+0*3.6*sin(\t r)},{0*3.6*cos(\t r)+1*3.6*sin(\t r)});
\draw [shift={(17.69,4.52)}] plot[domain=3.53:4.31,variable=\t]({1*3.99*cos(\t r)+0*3.99*sin(\t r)},{0*3.99*cos(\t r)+1*3.99*sin(\t r)});
\draw [shift={(15.49,-0.64)}] plot[domain=1.96:2.74,variable=\t]({1*3.94*cos(\t r)+0*3.94*sin(\t r)},{0*3.94*cos(\t r)+1*3.94*sin(\t r)});
\draw [shift={(10.18,4.74)}] plot[domain=5.13:5.86,variable=\t]({1*4.22*cos(\t r)+0*4.22*sin(\t r)},{0*4.22*cos(\t r)+1*4.22*sin(\t r)});
\draw [shift={(18.5,1.5)}] plot[domain=2.82:3.46,variable=\t]({1*4.74*cos(\t r)+0*4.74*sin(\t r)},{0*4.74*cos(\t r)+1*4.74*sin(\t r)});
\draw [shift={(10.55,1.66)}] plot[domain=0.37:1.21,variable=\t]({1*3.7*cos(\t r)+0*3.7*sin(\t r)},{0*3.7*cos(\t r)+1*3.7*sin(\t r)});
\draw [shift={(15.37,6.47)}] plot[domain=3.51:4.34,variable=\t]({1*3.74*cos(\t r)+0*3.74*sin(\t r)},{0*3.74*cos(\t r)+1*3.74*sin(\t r)});
\draw [->] (15.46,2.73) -- (15.58,2.73);
\draw [->] (15.55,3.27) -- (15.43,3.27);
\draw [->] (13.76,1.53) -- (13.76,1.41);
\draw [->] (14.27,1.43) -- (14.27,1.55);
\draw [->] (12.68,3.87) -- (12.77,3.78);
\draw [->] (13.19,4.25) -- (13.1,4.34);
\draw [->] (13.73,4.56) -- (13.73,4.44);
\draw [->] (14.27,4.46) -- (14.27,4.56);
\draw [->] (14.86,4.32) -- (14.78,4.24);
\draw [->] (15.25,3.8) -- (15.33,3.88);
\draw [->] (14.79,1.78) -- (14.88,1.69);
\draw [->] (15.32,2.13) -- (15.25,2.2);
\draw [->] (12.77,2.2) -- (12.69,2.12);
\draw [->] (13.09,1.68) -- (13.18,1.77);
\draw [->] (10.66,4.38) -- (10.58,4.53);
\draw (7.75,3.48) node[anchor=north west] {$C_q$};
\draw (10.54,3.64) node[anchor=north west] {$v$};
\draw (7.4,-4.8) node[anchor=north west] {$C_q$};
\draw (11.2,-4.8) node[anchor=north west] {$u$};
\draw [shift={(14,2)}] plot[domain=4.31:5.12,variable=\t]({1*7.62*cos(\t r)+0*7.62*sin(\t r)},{0*7.62*cos(\t r)+1*7.62*sin(\t r)});
\draw [shift={(14,-1)}] plot[domain=4.07:5.36,variable=\t]({1*5*cos(\t r)+0*5*sin(\t r)},{0*5*cos(\t r)+1*5*sin(\t r)});
\draw [shift={(14,-2.5)}] plot[domain=3.84:5.59,variable=\t]({1*3.91*cos(\t r)+0*3.91*sin(\t r)},{0*3.91*cos(\t r)+1*3.91*sin(\t r)});
\draw [shift={(14,-3.5)}] plot[domain=3.61:5.82,variable=\t]({1*3.35*cos(\t r)+0*3.35*sin(\t r)},{0*3.35*cos(\t r)+1*3.35*sin(\t r)});
\draw [shift={(14,-4)}] plot[domain=3.46:5.96,variable=\t]({1*3.16*cos(\t r)+0*3.16*sin(\t r)},{0*3.16*cos(\t r)+1*3.16*sin(\t r)});
\draw [shift={(14,-4.5)}] plot[domain=3.31:6.12,variable=\t]({1*3.04*cos(\t r)+0*3.04*sin(\t r)},{0*3.04*cos(\t r)+1*3.04*sin(\t r)});
\draw [shift={(14,-5)}] plot[domain=-3.14:0,variable=\t]({1*3*cos(\t r)+0*3*sin(\t r)},{0*3*cos(\t r)+1*3*sin(\t r)});
\draw [shift={(14,-5)}] plot[domain=0:3.14,variable=\t]({1*3.01*cos(\t r)+0*3.01*sin(\t r)},{0*3.01*cos(\t r)+1*3.01*sin(\t r)});
\draw [shift={(14,-5.5)}] plot[domain=0.16:2.98,variable=\t]({1*3.05*cos(\t r)+0*3.05*sin(\t r)},{0*3.05*cos(\t r)+1*3.05*sin(\t r)});
\draw [shift={(14,-6)}] plot[domain=0.32:2.82,variable=\t]({1*3.17*cos(\t r)+0*3.17*sin(\t r)},{0*3.17*cos(\t r)+1*3.17*sin(\t r)});
\draw [shift={(14,-6.5)}] plot[domain=0.46:2.68,variable=\t]({1*3.36*cos(\t r)+0*3.36*sin(\t r)},{0*3.36*cos(\t r)+1*3.36*sin(\t r)});
\draw [shift={(14,-7.54)}] plot[domain=0.7:2.44,variable=\t]({1*3.94*cos(\t r)+0*3.94*sin(\t r)},{0*3.94*cos(\t r)+1*3.94*sin(\t r)});
\draw [shift={(14,-9)}] plot[domain=0.93:2.21,variable=\t]({1*5.01*cos(\t r)+0*5.01*sin(\t r)},{0*5.01*cos(\t r)+1*5.01*sin(\t r)});
\draw [shift={(14,-12)}] plot[domain=1.16:1.98,variable=\t]({1*7.62*cos(\t r)+0*7.62*sin(\t r)},{0*7.62*cos(\t r)+1*7.62*sin(\t r)});
\draw [->] (13.83,-1.99) -- (14,-1.99);
\draw [->] (13.85,-2.45) -- (14,-2.45);
\draw [->] (13.84,-2.84) -- (14,-2.83);
\draw [->] (13.83,-3.14) -- (14,-3.14);
\draw [->] (13.85,-3.6) -- (14,-3.6);
\draw [->] (13.83,-4) -- (14,-3.99);
\draw [->] (13.84,-4.38) -- (14.01,-4.38);
\draw [->] (14.17,-5.61) -- (14.01,-5.62);
\draw [->] (14.17,-6) -- (14,-6);
\draw [->] (14.18,-6.4) -- (14,-6.41);
\draw [->] (14.18,-6.85) -- (14,-6.85);
\draw [->] (14.18,-7.16) -- (14,-7.16);
\draw [->] (14.18,-7.54) -- (14,-7.54);
\draw [->] (14.19,-7.99) -- (14,-8);
\draw [->] (16.03,-2.77) -- (16.13,-2.87);
\draw [->] (15.76,-3) -- (15.89,-3.11);
\draw [->] (15.55,-3.23) -- (15.69,-3.31);
\draw [->] (15.32,-3.4) -- (15.51,-3.49);
\draw [->] (15.05,-3.74) -- (15.21,-3.79);
\draw [->] (14.74,-4.05) -- (14.92,-4.08);
\draw [->] (14.49,-4.4) -- (14.6,-4.4);
\draw [->] (13.56,-5.6) -- (13.41,-5.59);
\draw [->] (13.26,-5.94) -- (13.08,-5.92);
\draw [->] (12.91,-6.25) -- (12.79,-6.21);
\draw [->] (12.67,-6.58) -- (12.5,-6.5);
\draw [->] (12.44,-6.75) -- (12.32,-6.68);
\draw [->] (12.23,-6.98) -- (12.11,-6.89);
\draw [->] (11.99,-7.23) -- (11.88,-7.12);
\draw [shift={(12.5,-0.5)}] plot[domain=1.17:1.98,variable=\t]({1*3.81*cos(\t r)+0*3.81*sin(\t r)},{0*3.81*cos(\t r)+1*3.81*sin(\t r)});
\draw [shift={(12.5,6.5)}] plot[domain=4.31:5.12,variable=\t]({1*3.81*cos(\t r)+0*3.81*sin(\t r)},{0*3.81*cos(\t r)+1*3.81*sin(\t r)});
\draw [->] (12.57,3.31) -- (12.42,3.31);
\draw [->] (12.44,2.69) -- (12.59,2.69);
\draw [shift={(10.5,-5.5)}] plot[domain=0.79:2.36,variable=\t]({1*0.71*cos(\t r)+0*0.71*sin(\t r)},{0*0.71*cos(\t r)+1*0.71*sin(\t r)});
\draw [shift={(10.5,-4.5)}] plot[domain=3.93:5.5,variable=\t]({1*0.71*cos(\t r)+0*0.71*sin(\t r)},{0*0.71*cos(\t r)+1*0.71*sin(\t r)});
\draw [->] (10.57,-4.8) -- (10.43,-4.8);
\draw [->] (10.44,-5.2) -- (10.58,-5.2);
\draw [dash pattern=on 5pt off 5pt] (7.56,-4.99) circle (2.44cm);
\draw [->] (9.46,-3.45) -- (9.26,-3.24);
\begin{scriptsize}
\fill [color=black] (11,3) circle (1.5pt);
\fill [color=black] (14,3) circle (3.0pt);
\fill [color=black] (14,6) circle (1.5pt);
\fill [color=black] (14,0) circle (1.5pt);
\fill [color=black] (17,3) circle (1.5pt);
\fill [color=black] (16.13,5.11) circle (1.5pt);
\fill [color=black] (16.12,0.86) circle (1.5pt);
\fill [color=black] (11.88,0.88) circle (1.5pt);
\fill [color=black] (11.88,5.12) circle (1.5pt);
\fill [color=black] (11,-5) circle (2.5pt);
\fill [color=black] (17.01,-5) circle (2.5pt);
\fill [color=black] (14.01,-5.62) circle (1.5pt);
\fill [color=black] (14,-6) circle (1.5pt);
\fill [color=black] (14,-6.41) circle (1.5pt);
\fill [color=black] (14,-6.85) circle (1.5pt);
\fill [color=black] (14,-7.16) circle (1.5pt);
\fill [color=black] (14,-8) circle (1.5pt);
\fill [color=black] (14,-7.54) circle (1.5pt);
\fill [color=black] (14.01,-4.38) circle (1.5pt);
\fill [color=black] (14,-3.99) circle (1.5pt);
\fill [color=black] (14,-3.6) circle (1.5pt);
\fill [color=black] (14,-3.14) circle (1.5pt);
\fill [color=black] (14,-2.83) circle (1.5pt);
\fill [color=black] (14,-2.45) circle (1.5pt);
\fill [color=black] (14,-1.99) circle (1.5pt);
\fill [color=black] (10,-5) circle (1.5pt);
\end{scriptsize}
\end{tikzpicture}
\caption{Example 1 of an MSD in which there is a vertex with high degree (input or output) and is not contained in the cycle $C_q$.}
\label{fig_contraejemplo1}
\end{figure}

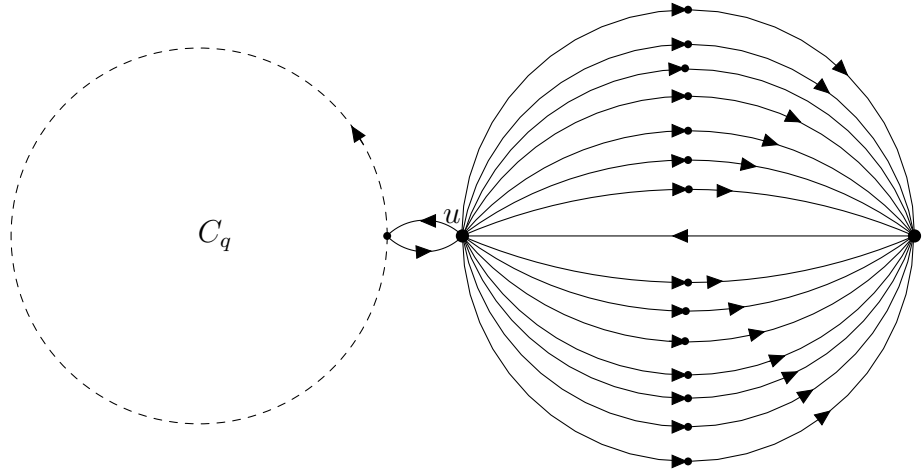
\begin{figure}[hp]
\centering
\begin{tikzpicture}[line cap=round,line join=round,>=triangle 45,x=1.0cm,y=1.0cm]
\clip(4.8,-15.17) rectangle (17.29,-8.84);
\draw [dash pattern=on 3pt off 3pt] (8.03,-4.97) circle (3cm);
\draw (7.62,-4.54) node[anchor=north west] {$C_q$};
\draw (10.5,-4.5) node[anchor=north west] {$v$};
\draw [shift={(14,2)}] plot[domain=4.31:5.12,variable=\t]({1*7.62*cos(\t r)+0*7.62*sin(\t r)},{0*7.62*cos(\t r)+1*7.62*sin(\t r)});
\draw [shift={(14,-1)}] plot[domain=4.07:5.36,variable=\t]({1*5*cos(\t r)+0*5*sin(\t r)},{0*5*cos(\t r)+1*5*sin(\t r)});
\draw [shift={(14,-2.5)}] plot[domain=3.84:5.59,variable=\t]({1*3.91*cos(\t r)+0*3.91*sin(\t r)},{0*3.91*cos(\t r)+1*3.91*sin(\t r)});
\draw [shift={(14,-3.5)}] plot[domain=3.61:5.82,variable=\t]({1*3.35*cos(\t r)+0*3.35*sin(\t r)},{0*3.35*cos(\t r)+1*3.35*sin(\t r)});
\draw [shift={(14,-4)}] plot[domain=3.46:5.96,variable=\t]({1*3.16*cos(\t r)+0*3.16*sin(\t r)},{0*3.16*cos(\t r)+1*3.16*sin(\t r)});
\draw [shift={(14,-4.5)}] plot[domain=3.31:6.12,variable=\t]({1*3.04*cos(\t r)+0*3.04*sin(\t r)},{0*3.04*cos(\t r)+1*3.04*sin(\t r)});
\draw [shift={(14,-5)}] plot[domain=-3.14:0,variable=\t]({1*3*cos(\t r)+0*3*sin(\t r)},{0*3*cos(\t r)+1*3*sin(\t r)});
\draw [shift={(14,-5)}] plot[domain=0:3.14,variable=\t]({1*3.01*cos(\t r)+0*3.01*sin(\t r)},{0*3.01*cos(\t r)+1*3.01*sin(\t r)});
\draw [shift={(14,-5.5)}] plot[domain=0.16:2.98,variable=\t]({1*3.05*cos(\t r)+0*3.05*sin(\t r)},{0*3.05*cos(\t r)+1*3.05*sin(\t r)});
\draw [shift={(14,-6)}] plot[domain=0.32:2.82,variable=\t]({1*3.17*cos(\t r)+0*3.17*sin(\t r)},{0*3.17*cos(\t r)+1*3.17*sin(\t r)});
\draw [shift={(14,-6.5)}] plot[domain=0.46:2.68,variable=\t]({1*3.36*cos(\t r)+0*3.36*sin(\t r)},{0*3.36*cos(\t r)+1*3.36*sin(\t r)});
\draw [shift={(14,-7.54)}] plot[domain=0.7:2.44,variable=\t]({1*3.94*cos(\t r)+0*3.94*sin(\t r)},{0*3.94*cos(\t r)+1*3.94*sin(\t r)});
\draw [shift={(14,-9)}] plot[domain=0.93:2.21,variable=\t]({1*5.01*cos(\t r)+0*5.01*sin(\t r)},{0*5.01*cos(\t r)+1*5.01*sin(\t r)});
\draw [shift={(14,-12)}] plot[domain=1.16:1.98,variable=\t]({1*7.62*cos(\t r)+0*7.62*sin(\t r)},{0*7.62*cos(\t r)+1*7.62*sin(\t r)});
\draw [->] (13.83,-1.99) -- (14,-1.99);
\draw [->] (13.85,-2.45) -- (14,-2.45);
\draw [->] (13.84,-2.84) -- (14,-2.83);
\draw [->] (13.83,-3.14) -- (14,-3.14);
\draw [->] (13.85,-3.6) -- (14,-3.6);
\draw [->] (13.83,-4) -- (14,-3.99);
\draw [->] (13.84,-4.38) -- (14.01,-4.38);
\draw [->] (14.17,-5.61) -- (14.01,-5.62);
\draw [->] (14.17,-6) -- (14,-6);
\draw [->] (14.18,-6.4) -- (14,-6.41);
\draw [->] (14.18,-6.85) -- (14,-6.85);
\draw [->] (14.18,-7.16) -- (14,-7.16);
\draw [->] (14.18,-7.54) -- (14,-7.54);
\draw [->] (14.19,-7.99) -- (14,-8);
\draw [->] (16.03,-2.77) -- (16.13,-2.87);
\draw [->] (15.76,-3) -- (15.89,-3.11);
\draw [->] (15.55,-3.23) -- (15.69,-3.31);
\draw [->] (15.32,-3.4) -- (15.51,-3.49);
\draw [->] (15.05,-3.74) -- (15.21,-3.79);
\draw [->] (14.74,-4.05) -- (14.92,-4.08);
\draw [->] (14.49,-4.4) -- (14.6,-4.4);
\draw [->] (13.56,-5.6) -- (13.41,-5.59);
\draw [->] (13.26,-5.94) -- (13.08,-5.92);
\draw [->] (12.91,-6.25) -- (12.79,-6.21);
\draw [->] (12.67,-6.58) -- (12.5,-6.5);
\draw [->] (12.44,-6.75) -- (12.32,-6.68);
\draw [->] (12.23,-6.98) -- (12.11,-6.89);
\draw [->] (11.99,-7.23) -- (11.88,-7.12);
\draw [->] (10.76,-3.74) -- (10.65,-3.52);
\draw (7.34,-11.65) node[anchor=north west] {$C_q$};
\draw (10.6,-11.5) node[anchor=north west] {$u$};
\draw [shift={(14,-5)}] plot[domain=4.31:5.12,variable=\t]({1*7.62*cos(\t r)+0*7.62*sin(\t r)},{0*7.62*cos(\t r)+1*7.62*sin(\t r)});
\draw [shift={(14,-8)}] plot[domain=4.07:5.36,variable=\t]({1*5*cos(\t r)+0*5*sin(\t r)},{0*5*cos(\t r)+1*5*sin(\t r)});
\draw [shift={(14,-9.5)}] plot[domain=3.84:5.59,variable=\t]({1*3.91*cos(\t r)+0*3.91*sin(\t r)},{0*3.91*cos(\t r)+1*3.91*sin(\t r)});
\draw [shift={(14,-10.5)}] plot[domain=3.61:5.82,variable=\t]({1*3.35*cos(\t r)+0*3.35*sin(\t r)},{0*3.35*cos(\t r)+1*3.35*sin(\t r)});
\draw [shift={(14,-11)}] plot[domain=3.46:5.96,variable=\t]({1*3.16*cos(\t r)+0*3.16*sin(\t r)},{0*3.16*cos(\t r)+1*3.16*sin(\t r)});
\draw [shift={(14,-11.5)}] plot[domain=3.31:6.12,variable=\t]({1*3.04*cos(\t r)+0*3.04*sin(\t r)},{0*3.04*cos(\t r)+1*3.04*sin(\t r)});
\draw [shift={(14,-12)}] plot[domain=-3.14:0,variable=\t]({1*3*cos(\t r)+0*3*sin(\t r)},{0*3*cos(\t r)+1*3*sin(\t r)});
\draw [shift={(14,-12)}] plot[domain=0:3.14,variable=\t]({1*3.01*cos(\t r)+0*3.01*sin(\t r)},{0*3.01*cos(\t r)+1*3.01*sin(\t r)});
\draw [shift={(14,-12.5)}] plot[domain=0.16:2.98,variable=\t]({1*3.05*cos(\t r)+0*3.05*sin(\t r)},{0*3.05*cos(\t r)+1*3.05*sin(\t r)});
\draw [shift={(14,-13.5)}] plot[domain=0.46:2.68,variable=\t]({1*3.36*cos(\t r)+0*3.36*sin(\t r)},{0*3.36*cos(\t r)+1*3.36*sin(\t r)});
\draw [shift={(14,-14.54)}] plot[domain=0.7:2.44,variable=\t]({1*3.94*cos(\t r)+0*3.94*sin(\t r)},{0*3.94*cos(\t r)+1*3.94*sin(\t r)});
\draw [shift={(14,-16)}] plot[domain=0.93:2.21,variable=\t]({1*5.01*cos(\t r)+0*5.01*sin(\t r)},{0*5.01*cos(\t r)+1*5.01*sin(\t r)});
\draw [shift={(14,-19)}] plot[domain=1.16:1.98,variable=\t]({1*7.62*cos(\t r)+0*7.62*sin(\t r)},{0*7.62*cos(\t r)+1*7.62*sin(\t r)});
\draw [->] (13.83,-8.99) -- (14,-8.99);
\draw [->] (13.85,-9.45) -- (14,-9.45);
\draw [->] (13.83,-10.14) -- (14,-10.14);
\draw [->] (13.85,-10.6) -- (14,-10.6);
\draw [->] (13.83,-11) -- (14,-10.99);
\draw [->] (13.84,-11.38) -- (14.01,-11.38);
\draw [->] (16.03,-9.77) -- (16.13,-9.87);
\draw [->] (15.76,-10) -- (15.89,-10.11);
\draw [->] (15.32,-10.4) -- (15.51,-10.49);
\draw [->] (15.05,-10.74) -- (15.21,-10.79);
\draw [->] (14.74,-11.05) -- (14.92,-11.08);
\draw [->] (14.49,-11.4) -- (14.6,-11.4);
\draw [->] (13.95,-15) -- (14,-15);
\draw [->] (13.94,-14.54) -- (14,-14.54);
\draw [->] (13.97,-14.16) -- (14,-14.16);
\draw [->] (13.97,-13.85) -- (14,-13.85);
\draw [->] (13.97,-13) -- (14,-13);
\draw [->] (13.97,-12.62) -- (14,-12.62);
\draw [->] (14.33,-12.61) -- (14.47,-12.6);
\draw [->] (14.57,-12.97) -- (14.76,-12.94);
\draw [->] (14.82,-13.32) -- (15.03,-13.27);
\draw [->] (15.08,-13.67) -- (15.3,-13.59);
\draw [->] (15.28,-13.89) -- (15.47,-13.8);
\draw [->] (15.51,-14.14) -- (15.67,-14.04);
\draw [->] (15.74,-14.44) -- (15.9,-14.32);
\draw (17.01,-12)-- (11,-12);
\draw [->] (14,-12) -- (13.78,-12);
\draw [shift={(10.5,-12.5)}] plot[domain=0.79:2.36,variable=\t]({1*0.71*cos(\t r)+0*0.71*sin(\t r)},{0*0.71*cos(\t r)+1*0.71*sin(\t r)});
\draw [shift={(10.5,-11.5)}] plot[domain=3.93:5.5,variable=\t]({1*0.71*cos(\t r)+0*0.71*sin(\t r)},{0*0.71*cos(\t r)+1*0.71*sin(\t r)});
\draw [dash pattern=on 3pt off 3pt] (7.5,-12) circle (2.5cm);
\draw [->] (10.59,-11.8) -- (10.43,-11.8);
\draw [->] (10.45,-12.21) -- (10.6,-12.2);
\draw [->] (9.68,-10.78) -- (9.52,-10.53);
\draw [shift={(14.01,-12.92)}] plot[domain=0.3:2.85,variable=\t]({1*3.14*cos(\t r)+0*3.14*sin(\t r)},{0*3.14*cos(\t r)+1*3.14*sin(\t r)});
\draw [->] (13.78,-9.78) -- (13.96,-9.77);
\draw [->] (13.75,-13.4) -- (14,-13.41);
\begin{scriptsize}
\fill [color=black] (11,-5) circle (2.5pt);
\fill [color=black] (17.01,-5) circle (2.5pt);
\fill [color=black] (14.01,-5.62) circle (1.5pt);
\fill [color=black] (14,-6) circle (1.5pt);
\fill [color=black] (14,-6.41) circle (1.5pt);
\fill [color=black] (14,-6.85) circle (1.5pt);
\fill [color=black] (14,-7.16) circle (1.5pt);
\fill [color=black] (14,-8) circle (1.5pt);
\fill [color=black] (14,-7.54) circle (1.5pt);
\fill [color=black] (14.01,-4.38) circle (1.5pt);
\fill [color=black] (14,-3.99) circle (1.5pt);
\fill [color=black] (14,-3.6) circle (1.5pt);
\fill [color=black] (14,-3.14) circle (1.5pt);
\fill [color=black] (14,-2.83) circle (1.5pt);
\fill [color=black] (14,-2.45) circle (1.5pt);
\fill [color=black] (14,-1.99) circle (1.5pt);
\fill [color=black] (11.03,-5) circle (1.5pt);
\fill [color=black] (11,-12) circle (2.5pt);
\fill [color=black] (17.01,-12) circle (2.5pt);
\fill [color=black] (14.01,-11.38) circle (1.5pt);
\fill [color=black] (14,-10.99) circle (1.5pt);
\fill [color=black] (14,-10.6) circle (1.5pt);
\fill [color=black] (14,-10.14) circle (1.5pt);
\fill [color=black] (14,-9.45) circle (1.5pt);
\fill [color=black] (14,-8.99) circle (1.5pt);
\fill [color=black] (14,-12.62) circle (1.5pt);
\fill [color=black] (14,-13.85) circle (1.5pt);
\fill [color=black] (14,-14.16) circle (1.5pt);
\fill [color=black] (14,-14.54) circle (1.5pt);
\fill [color=black] (14,-15) circle (1.5pt);
\fill [color=black] (10,-12) circle (1.5pt);
\fill [color=black] (13.96,-9.77) circle (1.5pt);
\fill [color=black] (13.97,-13) circle (1.5pt);
\fill [color=black] (13.96,-13.4) circle (1.5pt);
\end{scriptsize}
\end{tikzpicture}
\caption{Example 2 of an MSD in which there is a vertex with high degree (input or output) and is not contained in the cycle $C_q$.}
\label{fig_contraejemplo2}
\end{figure}

\begin{figure}[hp]
\centering
\begin{tikzpicture}[line cap=round,line join=round,>=triangle 45,x=1.0cm,y=1.0cm]
\clip(4.26,-0.61) rectangle (17.3,6.69);
\draw [dash pattern=on 6pt off 6pt] (8,3) circle (3cm);
\draw [shift={(18,4.5)}] plot[domain=2.78:3.5,variable=\t]({1*4.27*cos(\t r)+0*4.27*sin(\t r)},{0*4.27*cos(\t r)+1*4.27*sin(\t r)});
\draw [shift={(10,4.5)}] plot[domain=-0.36:0.36,variable=\t]({1*4.27*cos(\t r)+0*4.27*sin(\t r)},{0*4.27*cos(\t r)+1*4.27*sin(\t r)});
\draw [shift={(15.5,7)}] plot[domain=4.35:5.07,variable=\t]({1*4.27*cos(\t r)+0*4.27*sin(\t r)},{0*4.27*cos(\t r)+1*4.27*sin(\t r)});
\draw [shift={(15.5,-1)}] plot[domain=1.21:1.93,variable=\t]({1*4.27*cos(\t r)+0*4.27*sin(\t r)},{0*4.27*cos(\t r)+1*4.27*sin(\t r)});
\draw [shift={(10,1.5)}] plot[domain=-0.36:0.36,variable=\t]({1*4.27*cos(\t r)+0*4.27*sin(\t r)},{0*4.27*cos(\t r)+1*4.27*sin(\t r)});
\draw [shift={(17.31,1.79)}] plot[domain=1.91:2.79,variable=\t]({1*3.52*cos(\t r)+0*3.52*sin(\t r)},{0*3.52*cos(\t r)+1*3.52*sin(\t r)});
\draw [shift={(12.62,6.54)}] plot[domain=5.08:5.9,variable=\t]({1*3.8*cos(\t r)+0*3.8*sin(\t r)},{0*3.8*cos(\t r)+1*3.8*sin(\t r)});
\draw [shift={(12.74,-0.38)}] plot[domain=0.35:1.21,variable=\t]({1*3.6*cos(\t r)+0*3.6*sin(\t r)},{0*3.6*cos(\t r)+1*3.6*sin(\t r)});
\draw [shift={(17.69,4.52)}] plot[domain=3.53:4.31,variable=\t]({1*3.99*cos(\t r)+0*3.99*sin(\t r)},{0*3.99*cos(\t r)+1*3.99*sin(\t r)});
\draw [shift={(15.49,-0.64)}] plot[domain=1.96:2.74,variable=\t]({1*3.94*cos(\t r)+0*3.94*sin(\t r)},{0*3.94*cos(\t r)+1*3.94*sin(\t r)});
\draw [shift={(10.18,4.74)}] plot[domain=5.13:5.86,variable=\t]({1*4.22*cos(\t r)+0*4.22*sin(\t r)},{0*4.22*cos(\t r)+1*4.22*sin(\t r)});
\draw [shift={(18.5,1.5)}] plot[domain=2.82:3.46,variable=\t]({1*4.74*cos(\t r)+0*4.74*sin(\t r)},{0*4.74*cos(\t r)+1*4.74*sin(\t r)});
\draw [shift={(10.55,1.66)}] plot[domain=0.37:1.21,variable=\t]({1*3.7*cos(\t r)+0*3.7*sin(\t r)},{0*3.7*cos(\t r)+1*3.7*sin(\t r)});
\draw [shift={(15.37,6.47)}] plot[domain=3.51:4.34,variable=\t]({1*3.74*cos(\t r)+0*3.74*sin(\t r)},{0*3.74*cos(\t r)+1*3.74*sin(\t r)});
\draw [->] (15.46,2.73) -- (15.58,2.73);
\draw [->] (15.55,3.27) -- (15.43,3.27);
\draw [->] (13.76,1.53) -- (13.76,1.41);
\draw [->] (14.27,1.43) -- (14.27,1.55);
\draw [->] (12.68,3.87) -- (12.77,3.78);
\draw [->] (13.19,4.25) -- (13.1,4.34);
\draw [->] (13.73,4.56) -- (13.73,4.44);
\draw [->] (14.27,4.46) -- (14.27,4.56);
\draw [->] (14.86,4.32) -- (14.78,4.24);
\draw [->] (15.25,3.8) -- (15.33,3.88);
\draw [->] (14.79,1.78) -- (14.88,1.69);
\draw [->] (15.32,2.13) -- (15.25,2.2);
\draw [->] (12.77,2.2) -- (12.69,2.12);
\draw [->] (13.09,1.68) -- (13.18,1.77);
\draw [->] (10.66,4.38) -- (10.58,4.53);
\draw (7.75,3.2) node[anchor=north west] {$C_q$};
\draw (13.3,3.2) node[anchor=north west] {$u$};
\draw [dash pattern=on 6pt off 6pt] (4.99,-4.98) circle (3cm);
\draw (4.67,-4.36) node[anchor=north west] {$C_q$};
\draw (7.22,-4.21) node[anchor=north west] {$v$};
\draw [shift={(14,2)}] plot[domain=4.31:5.12,variable=\t]({1*7.62*cos(\t r)+0*7.62*sin(\t r)},{0*7.62*cos(\t r)+1*7.62*sin(\t r)});
\draw [shift={(14,-1)}] plot[domain=4.07:5.36,variable=\t]({1*5*cos(\t r)+0*5*sin(\t r)},{0*5*cos(\t r)+1*5*sin(\t r)});
\draw [shift={(14,-2.5)}] plot[domain=3.84:5.59,variable=\t]({1*3.91*cos(\t r)+0*3.91*sin(\t r)},{0*3.91*cos(\t r)+1*3.91*sin(\t r)});
\draw [shift={(14,-3.5)}] plot[domain=3.61:5.82,variable=\t]({1*3.35*cos(\t r)+0*3.35*sin(\t r)},{0*3.35*cos(\t r)+1*3.35*sin(\t r)});
\draw [shift={(14,-4)}] plot[domain=3.46:5.96,variable=\t]({1*3.16*cos(\t r)+0*3.16*sin(\t r)},{0*3.16*cos(\t r)+1*3.16*sin(\t r)});
\draw [shift={(14,-4.5)}] plot[domain=3.31:6.12,variable=\t]({1*3.04*cos(\t r)+0*3.04*sin(\t r)},{0*3.04*cos(\t r)+1*3.04*sin(\t r)});
\draw [shift={(14,-5)}] plot[domain=-3.14:0,variable=\t]({1*3*cos(\t r)+0*3*sin(\t r)},{0*3*cos(\t r)+1*3*sin(\t r)});
\draw [shift={(14,-5)}] plot[domain=0:3.14,variable=\t]({1*3.01*cos(\t r)+0*3.01*sin(\t r)},{0*3.01*cos(\t r)+1*3.01*sin(\t r)});
\draw [shift={(14,-5.5)}] plot[domain=0.16:2.98,variable=\t]({1*3.05*cos(\t r)+0*3.05*sin(\t r)},{0*3.05*cos(\t r)+1*3.05*sin(\t r)});
\draw [shift={(14,-6)}] plot[domain=0.32:2.82,variable=\t]({1*3.17*cos(\t r)+0*3.17*sin(\t r)},{0*3.17*cos(\t r)+1*3.17*sin(\t r)});
\draw [shift={(14,-6.5)}] plot[domain=0.46:2.68,variable=\t]({1*3.36*cos(\t r)+0*3.36*sin(\t r)},{0*3.36*cos(\t r)+1*3.36*sin(\t r)});
\draw [shift={(14,-7.54)}] plot[domain=0.7:2.44,variable=\t]({1*3.94*cos(\t r)+0*3.94*sin(\t r)},{0*3.94*cos(\t r)+1*3.94*sin(\t r)});
\draw [shift={(14,-9)}] plot[domain=0.93:2.21,variable=\t]({1*5.01*cos(\t r)+0*5.01*sin(\t r)},{0*5.01*cos(\t r)+1*5.01*sin(\t r)});
\draw [shift={(14,-12)}] plot[domain=1.16:1.98,variable=\t]({1*7.62*cos(\t r)+0*7.62*sin(\t r)},{0*7.62*cos(\t r)+1*7.62*sin(\t r)});
\draw [->] (13.83,-1.99) -- (14,-1.99);
\draw [->] (13.85,-2.45) -- (14,-2.45);
\draw [->] (13.84,-2.84) -- (14,-2.83);
\draw [->] (13.83,-3.14) -- (14,-3.14);
\draw [->] (13.85,-3.6) -- (14,-3.6);
\draw [->] (13.83,-4) -- (14,-3.99);
\draw [->] (13.84,-4.38) -- (14.01,-4.38);
\draw [->] (14.17,-5.61) -- (14.01,-5.62);
\draw [->] (14.17,-6) -- (14,-6);
\draw [->] (14.18,-6.4) -- (14,-6.41);
\draw [->] (14.18,-6.85) -- (14,-6.85);
\draw [->] (14.18,-7.16) -- (14,-7.16);
\draw [->] (14.18,-7.54) -- (14,-7.54);
\draw [->] (14.19,-7.99) -- (14,-8);
\draw [->] (16.03,-2.77) -- (16.13,-2.87);
\draw [->] (15.76,-3) -- (15.89,-3.11);
\draw [->] (15.55,-3.23) -- (15.69,-3.31);
\draw [->] (15.32,-3.4) -- (15.51,-3.49);
\draw [->] (15.05,-3.74) -- (15.21,-3.79);
\draw [->] (14.74,-4.05) -- (14.92,-4.08);
\draw [->] (14.49,-4.4) -- (14.6,-4.4);
\draw [->] (13.56,-5.6) -- (13.41,-5.59);
\draw [->] (13.26,-5.94) -- (13.08,-5.92);
\draw [->] (12.91,-6.25) -- (12.79,-6.21);
\draw [->] (12.67,-6.58) -- (12.5,-6.5);
\draw [->] (12.44,-6.75) -- (12.32,-6.68);
\draw [->] (12.23,-6.98) -- (12.11,-6.89);
\draw [->] (11.99,-7.23) -- (11.88,-7.12);
\draw [->] (7.72,-3.74) -- (7.61,-3.52);
\draw [shift={(9.5,-8)}] plot[domain=1.11:2.04,variable=\t]({1*3.35*cos(\t r)+0*3.35*sin(\t r)},{0*3.35*cos(\t r)+1*3.35*sin(\t r)});
\draw [shift={(9.5,-2)}] plot[domain=4.24:5.18,variable=\t]({1*3.36*cos(\t r)+0*3.36*sin(\t r)},{0*3.36*cos(\t r)+1*3.36*sin(\t r)});
\draw [->] (9.46,-4.65) -- (9.6,-4.65);
\draw [->] (9.54,-5.36) -- (9.38,-5.36);
\draw [shift={(12.5,-0.5)}] plot[domain=1.17:1.98,variable=\t]({1*3.81*cos(\t r)+0*3.81*sin(\t r)},{0*3.81*cos(\t r)+1*3.81*sin(\t r)});
\draw [shift={(12.5,6.5)}] plot[domain=4.31:5.12,variable=\t]({1*3.81*cos(\t r)+0*3.81*sin(\t r)},{0*3.81*cos(\t r)+1*3.81*sin(\t r)});
\draw [->] (12.57,3.31) -- (12.42,3.31);
\draw [->] (12.44,2.69) -- (12.59,2.69);
\begin{scriptsize}
\fill [color=black] (11,3) circle (1.5pt);
\fill [color=black] (14,3) circle (3.0pt);
\fill [color=black] (14,6) circle (1.5pt);
\fill [color=black] (14,0) circle (1.5pt);
\fill [color=black] (17,3) circle (1.5pt);
\fill [color=black] (16.13,5.11) circle (1.5pt);
\fill [color=black] (16.12,0.86) circle (1.5pt);
\fill [color=black] (11.88,0.88) circle (1.5pt);
\fill [color=black] (11.88,5.12) circle (1.5pt);
\fill [color=black] (11,-5) circle (2.5pt);
\fill [color=black] (17.01,-5) circle (2.5pt);
\fill [color=black] (14.01,-5.62) circle (1.5pt);
\fill [color=black] (14,-6) circle (1.5pt);
\fill [color=black] (14,-6.41) circle (1.5pt);
\fill [color=black] (14,-6.85) circle (1.5pt);
\fill [color=black] (14,-7.16) circle (1.5pt);
\fill [color=black] (14,-8) circle (1.5pt);
\fill [color=black] (14,-7.54) circle (1.5pt);
\fill [color=black] (14.01,-4.38) circle (1.5pt);
\fill [color=black] (14,-3.99) circle (1.5pt);
\fill [color=black] (14,-3.6) circle (1.5pt);
\fill [color=black] (14,-3.14) circle (1.5pt);
\fill [color=black] (14,-2.83) circle (1.5pt);
\fill [color=black] (14,-2.45) circle (1.5pt);
\fill [color=black] (14,-1.99) circle (1.5pt);
\fill [color=black] (7.99,-5) circle (1.5pt);
\end{scriptsize}
\end{tikzpicture}
\caption{Example 3 of an MSD in which there is a vertex with high degree (input or output) and is not contained in the cycle $C_q$.}
\label{fig_contraejemplo3}
\end{figure}
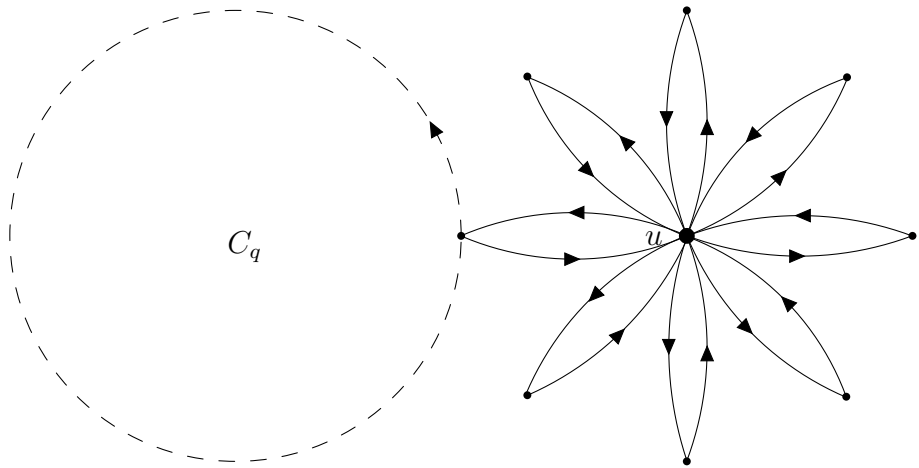

Corollary \ref{cor_lamda_maxgrado} is very useful when a vertex $u$ with high degree (input or output) is contained in the cycle $C_q$ (see examples in Figures \ref{fig_ejemplo1} and \ref{fig_ejemplo2}). However, if the vertex $u$ is not contained in the cycle, the number of linear vertices contained in the MSD could be much higher than the number of linear vertices obtained with this bound (see examples in Figures \ref{fig_contraejemplo1}, \ref{fig_contraejemplo2} and \ref{fig_contraejemplo3}).

%%%%%%%%%%%%%%%%%%%%%%%%%%%%%%%%%%%%%%%%%%%%%%%%%%%%%%%%%%%%%%%%%%%%%%%%%%%%%%%%%%%%%%%%%%%%%%%
%**************************************************************************************
\pagebreak
\color{black}
\begin{pro}\label{prop_1}
Let $D=(V,A)$ be an MSD and let $C_q$ be a cycle of length $q$ contained in $D$. Then $q \leq 2n-m$.
\end{pro}
\begin{pf}

\color{black} We obtain an MSD $D'$ by contraction of $C_q$ in an unique vertex $v'$, and then $n'=n-q+1$ and $m'=m-q$. Hence, since
$$m' \leq 2(n'-1),$$
we obtain
$$m-q \leq 2(n-q),$$
and finally
$$q \leq 2n-m.$$
The proof is completed.
\end{pf}

%**************************************************************************************
\subsection{Other properties of MSDs} \label{subsection31}

In \cite{BJG,cita2}, some results about ear decomposition are proved. We have used these previous results to show the next properties of MSDs.

\begin{df} \label{dfchain}
Let $D=(V,A)$ be an MSD of order $n \geq 2$ and let $v_1,\dots,v_l$ be a path contained in $D$. We say that the $v_1v_l$-path is a \emph{chain} with length $l$ if $d^-(v_i)=d^+(v_i)=1$ for all $1 \leq i \leq l$.
\end{df}

Note that an isolated linear vertex is a chain of length 1.

\begin{df} \label{dfextchain}
Let $D=(V,A)$ be an MSD of order $n \geq 2$, let $v_1,\dots,v_l$ be a chain contained in $D$ and let $D'$ be the digraph obtained from $D$ by the elimination of the $v_1v_l$-path. We say that the $v_1v_l$-path is an \emph{external chain} with length $l$ if $D'$ preserves the strongly connection.
\end{df}

%2.3 a Lemma 8%%%%%%%%%%%%%%%%%%%%%%%%%%%%%%%%%%%%%%%%%%%%%%%%%%%%%%%%%%%%%%%

\begin{lem}\label{lema2_3}
Let $D=(V,A)$ be an MSD of order $n \geq 2$, $C_q$ a cycle contained in $D$ such that $D \neq C_q$. Then, in $D$ there exists at least one external chain.
\end{lem}
\begin{pf}
We use the ear decomposition showed in ~\cite[Theorem 20]{cita2}, in a similar way as it was used in the proof of \color{black}the property that affirms that an MSD factors into a rooted spanning tree and a forest of reversed rooted trees\color{black}.

Let us consider an ear decomposition of $D$, $\mathcal{E}={P_0,\dots,P_k}$. Since $D$ is an MSD, each ear $P_j$ ($0 \leq j \leq k$) contains at least one new vertex and two new arcs, with respect to $\bigcup_{i=0}^{j-1}V_i$ and $\bigcup_{i=0}^{j-1}A_i$, respectively.

Then, it is clear that the last ear $P_k=v_{0}^{k} \dots v_{s_k}^{k}$ completes the construction of $D$, and $Q_k=v_1^k \dots v_{s_k-1}^k$ is a chain of linear vertices, whose first and last vertex are joined to vertices of a minimal and strongly connected digraph $D'$. Hence, $D'=D-Q_k$ is an MSD, and therefore $Q_k$ is an external chain of length $l = s_k-1 \geq 1$. Trivially we can say that if $D=Cq=P_0$, then there is no external chain contained in $D$.
The proof is completed.
\end{pf}

%2.2   2.4%%%%%%%%%%%%%%%%%%%%%%%%%%%%%%%%%%%%%%%%%%%%%%%%%%%%%%%%%%%%%%%%%%%%%%%%%
Note that $D'$ is an MSD with $n-l$ vertices and $m-l-1$ arcs. Note also that if $P_0=C_q$, with $C_q$ a maximal length cycle contained in $D$ and there exists any external chain with length $l \geq 1$ not contained in $C_q$, then $q \leq n-l$.

%2.1 a Lemma 9%%%%%%%%%%%%%%%%%%%%%%%%%%%%%%%%%%%%%%%%%%%%%%%%%%%%%%%%%%%%%%%
\begin{lem}\label{lema2_3}
Let $D=(V,A)$ be an MSD and let $v_1v_l$-path be a chain contained in $D$ with length $l < n$. Then, the contraction of all vertices of the $v_1v_l$-path in a unique vertex preserves the minimality, that is, it produces another MSD $D'$ with $n-l+1$ vertices and $m-l+1$ arcs.
\end{lem}
\begin{pf}
Let $D'$ be the digraph obtained by the contraction of all vertices of the $v_1v_l$-path in a unique vertex $v'$. Let $n'$ be the number of vertices and $m'$ the number of arcs of $D'$. In $D'$ all vertices of $v_1v_l$-path are suppressed but it contains the vertex $v' \notin D$, and then $n'=n-l+1$. Since $d^-(v_i)=d^+(v_i)=1$ for all $1\leq i\leq l$, we have $m'=m-l+1$. Now, let us assume that there are transitive arcs in $D'$. If we expand $v'$, these transitive arcs would also exist in $D$, contradicting the minimality of $D$. Hence, $D'$ is minimal. Since $n>l$, then a vertex $w \notin v_1v_l$-path, exists also in $D'$, and $D'$ contains a $wv'$-path and a $v'w$-path. Therefore $D'$ is strongly connected. The proof is completed.
\end{pf}

%2.5%%%%%%%%%%%%%%%%%%%%%%%%%%%%%%%%%%%%%%%%%%%%%%%%%%%%%%%%%%%%%%%
\begin{lem}\label{lema2_3}
Let $D=(V,A)$ be an MSD such that $D$ is not a cycle. Then, there is not a cycle in $D$ that contains all linear vertices of $D$.
\end{lem}
\begin{pf}
Let us suppose that $C_q$ contains all linear vertices of $D$. We can obtain an MSD $D'$ by contraction of $C_q$ in an unique vertex $v'$. We know that $D'$ must contain at least two linear vertices, and at least one of them is different from $v'$. Then, it is clear that there exists at least one linear vertex that is contained in $D$ but is not contained in $C_q$.
The proof is completed.
\end{pf}

%2.5%%%%%%%%%%%%%%%%%%%%%%%%%%%%%%%%%%%%%%%%%%%%%%%%%%%%%%%%%%%%%%%

Let $D$ be an MSD such that $D$ is not a cycle, and $\lambda$ be the number of linear vertices contained in $D$. From the lemma above, it is trivial to see that a cycle $C_q$ contained in $D$ will contain at most $\lambda - 1$ linear vertices of $D$.

%**************************************************************************************
\section{Upper bounds for the coefficients of the characteristical polynomial of MSDs} \label{section4}
%**************************************************************************************
In \cite{cita2}, some results about bounds of the coefficients of the characteristic polynomial of an MSD are proved. In particular, it is shown that the independent term must be $1$, $0$ or $-1$. We follow the lines of that proof to generalize that bound.

\begin{pro}
Let $D=(V,A)$ be an MSD, and let $x^n+k_1x^{n-1}+\cdots+k_i x^{n-1}+\cdots+k_{n-1}x+k_n$
be the characteristic polynomial of the adjacency matrix of $D$. Then
\[
|k_i|\leq \binom{n}{i}
\]
\end{pro}
\begin{pf}
We claim that any subset of $i$ vertices can be covered by disjoint cycles in at most one manner. In fact, take any subset $A\subset V$,  with $|A|=i$, consider be the subdigraph $D'$ generated by that $A$. Now, $D'$ is a subdigraph of an MSD, so it has no transitive arcs. If it is not strongly connected, we can add arcs, one by one, until we obtain a strongly connected digraph $D''$, that would be minimal. Therefore, $D''$ would be an MSD, and the aforementioned result of \cite{cita2} implies that there is at most one covering of the vertices of $D''$ (that is, of $A$) by disjoint cycles.

The coefficients theorem for digraphs (see, for instance, \cite[Theorem 2.1]{CDS}) allows us to finish the proof.
\end{pf}
%**************************************************************************************
\section{MSD properties associated to results of algorithms complexity} \label{section5}
It was well known that minimality is a very strict condition in the family of strong digraphs implying, for instance, the size limitation $n\leq m\leq 2(n-1)$. As we have seen in previous sections, MSDs also exhibit strong constraints on the number of linear vertices and maximum in- and outdegrees of vertices, regarding the length of the longest directed
cycle. Unfortunately, these constraints are not enough to construct an efficient algorithm finding the longest cycle in an MSD.

\begin{figure}[htb]
\centering
\begin{tikzpicture}[line cap=round,line join=round,>=triangle 45,x=1.0cm,y=1.0cm]
\clip(0.4,0.29) rectangle (6.64,4.44);
\draw [->] (1,1) -- (0.97,3.01);
\draw [->] (0.97,3.01) -- (3,3);
\draw [->] (3,3) -- (3,1);
\draw [->] (3,1) -- (1,1);
\draw [->] (3,1) -- (5,1);
\draw [->] (5,1) -- (5,3);
\draw [->] (5,3) -- (3,3);
\draw [->] (5,1) -- (6,4);
\draw [->] (6,4) -- (0.97,3.01);
\begin{scriptsize}
\fill [color=black] (1,1) circle (1.5pt);
\draw[color=black] (0.81,0.8) node {$u_7$};
\fill [color=black] (0.97,3.01) circle (1.5pt);
\draw[color=black] (0.81,3.17) node {$u_1$};
\fill [color=black] (3,3) circle (1.5pt);
\draw[color=black] (2.95,3.17) node {$u_2$};
\fill [color=black] (3,1) circle (1.5pt);
\draw[color=black] (2.95,0.8) node {$u_3$};
\fill [color=black] (5,1) circle (1.5pt);
\draw[color=black] (5,0.8) node {$u_4$};
\fill [color=black] (5,3) circle (1.5pt);
\draw[color=black] (5,3.17) node {$u_6$};
\fill [color=black] (6,4) circle (1.5pt);
\draw[color=black] (6.1,4.15) node {$u_5$};
\end{scriptsize}
\end{tikzpicture}
\caption{Example of an MSD that contains three external chains.}
\label{fig1}
\end{figure}
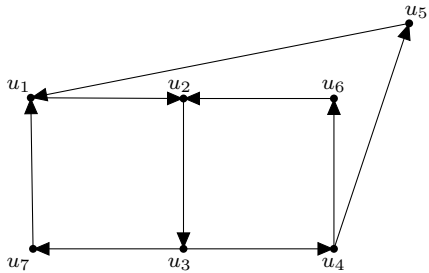

A proof that an MSD can be converted into a directed cycle by successively eliminating of external chains is given in \cite{cita4}. However, this process does not guarantee that the resulting directed cycle will have a maximum length. Figure \ref{fig1} shows an MSD where the longest directed cycle is given by $u_1,u_2,u_3,u_4,u_5,u_1$, but this cycle will be obtained only in the case that the external chains eliminated are those formed by $u_6$-path and $u_7$-path. Nevertheless, there is no an efficient algorithm that can determine the deletion of these chains and the non-deletion of the external chain formed by the vertex $u_5$-path, because if this chain is deleted then the longest cycle of the MSD will also have been eliminated.

\begin{thm}\label{complex1}
Computing a cycle with maximal length in an MSD is an NP-hard problem.
\end{thm}
\begin{pf}
We can reduce the problem of computing a cycle with maximal length in a strongly connected digraph to the problem of computing a cycle with maximal length in an MSD.

Let $D'=(V',A')$ be a strong digraph. We can build an MSD $D=(V,A)$ from $D'$ as follows. For each arc $v'_i v'_j\in A'$ we add an intermediate vertex $v_{ij}$. We thus obtain
\[
    V=V'\cup \{ v_{ij} \text{ } \mid \text{ } v'_iv'_j\in A'\},
\]
\[
    A=\{v'_iv_{ij} \text{ } \mid \text{ } v'_iv'_j\in A'\}\cup \{v_{ij}v'_j \text{ } \mid \text{ } v'_iv'_j\in A'\}.
\]
Note that the strong connection of $D'$ implies that $D$ is trivially strongly connected. Note also
that no arc of $D$ can be transitive, since every arc has a linear vertex $v_{ij}$ as start- or endpoint. Hence, $D$ is in fact an MSD.

Now, we remark that there is a one-to-one correspondence between cycles in $D$ and cycles in $D'$: for every cycle $C'_q$ in $D'$, a cycle $C_{2q}$ arises in $D$, and all the cycles in $D$ are generated in this way.

We conclude that any algorithm allowing us to compute the longest cycle of an MSD would then be able to compute the longest cycle of any SD, too. Since the problem of compute the longest cycle in a strongly connected digraph is NP-Hard \cite{BJG}, then the theorem is proved.
%%Añadir referencia de longest cycle in SD is NP-hard
\end{pf}

\begin{thm}\label{complex2}
Let $D=(V,A)$ be an MSD. Finding a cycle contained in $D$ with length $2n-m$ is an NP-complete problem.
\end{thm}
\begin{pf}
We can reduce the problem of determining if a digraph is Hamiltonian to the problem of determining if an MSD has a cycle of length $2n-m$.

Let $D'=(V',A')$ be a digraph. If $D'$ is strongly connected, the same procedure used in previous proof yields an MSD $D=(V,A)$ (if $D'$ is not strongly connected, then it can't be Hamiltonian). The order of $D$ verifies $n=n'+m'$, and the size holds $m=2m'$. Hence, finding a cycle in $D$ with length $2n-m=2(n'+m')-2m'=2n'$ would imply to find a $n'$-cycle in $D'$, that is, to determine if $D$ is Hamiltonian. Since to determine if a digraph is Hamiltonian is a NP-complete problem, the theorem is proved.
   %%añadir referencia de que decidir si un digrafo es hamiltoniano es NPC
   %%Thus, Degrees Theorem allow to proof Cycles Theorem (over lower bound of linear vertices).
\end{pf}

\section{Conclusions} \label{conclusiones}
In this work, we have found some new properties regarding MSDs. The first set of properties has to do with the number of linear vertices in an MSD. We have seen that the existence of a vertex with high in- or outdegree implies a high number of linear vertices. Furthermore, we have use this fact to give a simpler proof of the lower bound of linear vertices that we obtained in \cite{cita5}, where the existence of a $q$ cycle implied at least $\lfloor (q+1)/2 \rfloor$ linear vertices. We have also proved that chains of consecutive linear vertices in an MSD can be contracted, without loss of minimality. We feel that further research along these lines could give, from one side, a result linking maximal cycle lengths, maximal in- or outdegrees and improved estimations of number of linear vertices; and also a better understanding of cycle properties that could lead to spectral properties, such as characterization of polynomials that can be realized as characteristic polynomials of MSDs. In this regard, we have proved a bound for the coefficients of such polynomials, advancing along the lines given in \cite{cita2}.

Since the number of linear vertices in an MSD is easily computed, we wanted to explore the possibility that the maximal length of a cycle could be bounded so as to allow to construct a polynomial complexity algorithm to find the longest cycle. Unfortunately, that is not the case and we have proved that the search of a maximal length cycle in an MSD is NP-hard. Still, it could be interesting to look for a subset of MSDs for with the search for maximal length cycles can be performed efficiently. This kind of result could arise, also, by the further study of the properties that we have pointed out in the paragraph above.

%\section*{References}

%\bibliography{mybibfile}

\end{document}